\documentclass[final,3p,times]{elsarticle}

\usepackage{amssymb}
\usepackage{amsthm}
\usepackage{lipsum}
\usepackage{amsmath}
\usepackage{geometry}
\usepackage{graphicx}
\usepackage{epstopdf}
\usepackage{xcolor}
\usepackage{diagbox}
\usepackage{multirow}
\usepackage{booktabs}
\usepackage{amsopn}
\usepackage{longtable}
\usepackage{hyperref}
\usepackage{enumitem}
\hypersetup{hypertex=true,
            colorlinks=true,
            linkcolor=blue,
            anchorcolor=blue,
            citecolor=blue}
\newtheorem{thm}{Theorem}[section]
\newtheorem{lem}{Lemma}[section]
\newtheorem{prop}{Proposition}[section]
\newdefinition{rmk}{Remark}
\newproof{pf}{Proof}
\allowdisplaybreaks[4]
\usepackage{mathrsfs}
\numberwithin{figure}{section}
\newcommand{\rred}[1]{\textcolor{red}{#1}}

\begin{document}

\begin{frontmatter}



\title{The contour integral method for Feynman-Kac equation with two internal states}



\author[cor2]{Fugui Ma}
\author[cor3]{Lijing zhao}
\author[cor2]{Yejuan Wang}
\author[cor2]{Weihua Deng\corref{cor1}}
\ead{dengwh@lzu.edu.cn}
\cortext[cor1]{Corresponding author}
\address[cor2]{School of Mathematics and Statistics, Gansu Key Laboratory of Applied Mathematics and Complex Systems,
Lanzhou University, Lanzhou 730000, Peoples Republic of China.}
\address[cor3]{I. Research and Development Institute of Northwestern Polytechnical University in Shenzhen. II. School of Mathematics and Statistics, Northwestern Polytechnical University, Xian, 710129, Peoples Republic of China.}

\begin{abstract}
We develop the contour integral method for numerically solving the Feynman-Kac equation with two internal states [P. B. Xu and W. H. Deng, Math.
Model. Nat. Phenom., 13 (2018), 10], describing the functional distribution of particle's internal states. The striking benefits are obtained, including spectral accuracy, low computational complexity, small memory requirement, etc. We perform the error estimates and stability analyses, which are confirmed by numerical experiments.
\end{abstract}




\end{frontmatter}


\section{Introduction}
\label{sec1}
The weak singularity at starting point of the solution and the non-locality of the time evolution operator of the Feynman-Kac equation \cite{xu18} bring challenges on computational efficiency in numerically solving the equation. One of the most effective techniques to overcome the challenges is to analytically get the solution in frequency domain and then numerically do the inverse Laplace transform  \cite{Piessens75,Piessens76,Taiwo95}. The contour integral method (CIM) is an efficient numerical method for solving the inverse Laplace transform \cite{Talbot79,Taiwo95,Weidemantrefe07}.

Let us briefly introduce this method through the following toy model. Consider the time fractional initial value problem \cite{Podlubny99,Deng07JCAM}
\begin{equation}\label{eq:iniproblem1}
   ~_{0}^{C}D_{t}^{\alpha}\mathbf{u}(t) + A\mathbf{u}(t) = \mathbf{b}(t),\
                                 \mathbf{u}(0)=\mathbf{u}_{0},
\end{equation}
where $A$ is a matrix and $~_{0}^{C}D_{t}^{\alpha}$ is the Caputo fractional derivative \cite{Podlubny99} with $t\in (0,T]$ and $0<\alpha<1$.
Taking Laplace transform on (\ref{eq:iniproblem1}), one can get the solution of the system in Laplace space, namely,
\begin{equation}\label{eq:cost}
   \widehat{\mathbf{u}}(z)=(z^{\alpha}I+A)^{-1} \left(z^{\alpha-1}\mathbf{u}_{0}+\widehat{\mathbf{b}}(z)\right).
\end{equation}
Further by performing the inverse Laplace transform on the solution in Laplace space, one get the solution of the system (\ref{eq:iniproblem1}), i.e.,
\begin{equation}\label{eq:Lapalce}
  \mathbf{u}(t)=\frac{1}{2\pi i}\int_{\sigma-i\infty}^{\sigma+i\infty}e^{zt}\widehat{\mathbf{u}}(z)dz, \quad \sigma>\sigma_{0},
\end{equation}
where $\sigma_{0}$ is called convergent abscissa. In practice, due to the complexity of $\widehat{\mathbf{u}}(z)$ and the high dimension of matrix $A$, it is hard to get the analytical solution of (\ref{eq:iniproblem1}) by using the inverse Laplace transform (\ref{eq:Lapalce}). Hence, numerical methods are usually used to approximate (\ref{eq:Lapalce}).
The CIM is one of the most efficient numerical methods to solve this indefinite integral.

The earliest discussion for CIM seems appeared in \cite{Talbot79} by A. Talbot. Then, J. A. C. Weideman and other researchers gradually improved it and made it more efficient and widely applicable (see, e.g., \cite{Lee06,Sloan99}, etc). The basic idea of the CIM method is to deform the integral line, since  the original integration path of the inverse Laplace transform is a vertical line from negative infinity to positive infinity, which has many numerical challenges, e.g., the high frequency oscillation of the integrand.
Fortunately, by deforming the vertical line into a curve, which starts and ends in the left half complex plane,
the exponential decay of the integrand can be obtained by the exponential factor $e^{zt}$. Such a deformed line and the exponential decay of the integrand  make it possible for an abnormal integral to be solved numerically, and Cauchy integral theorem can ensure that such a deformation can be carried out. More specifically, after deforming the integral line of (\ref{eq:Lapalce}) into a contour,
which satisfies $Re(z)\rightarrow-\infty$ at each end, then, the exponential factor $e^{zt}$ forces a rapid decay of the integrand as $Re(z)\rightarrow-\infty$, which greatly benefits to the convergence speed of numerical integral methods for solving the inverse Laplace transform.

Based on the aforementioned idea, we return to the original time fractional initial value problem (\ref{eq:iniproblem1}). Suppose that there is an appropriate contour (for (\ref{eq:iniproblem1})) parameterized by
\begin{equation}\label{eq:contour}
\Gamma:z=z(\phi),\ -\infty<\phi<\infty.
\end{equation}
Then the solution of (\ref{eq:Lapalce}) can be rewritten as
\begin{equation}\label{eq:wcontour}
  \mathbf{u}(t)=\frac{1}{2\pi i}\int_{-\infty}^{\infty}e^{z(\phi)t}\widehat{\mathbf{u}}(z(\phi))z'(\phi)d\phi.
\end{equation}
Approximating it by the trapezoidal rule with uniform step-length $h$, there is
\begin{equation}\label{eq:approx}
 \mathbf{u}(t) \approx \frac{h}{2\pi i}\sum\limits_{k = -\infty}^{\infty}e^{z_{k}t}\widehat{\mathbf{u}}_{k}z_{k}',
\end{equation}
where $z_{k}=z(\phi_{k})$, $z'_{k}=z'(\phi_{k})$, $\widehat{\mathbf{u}}_{k}=\widehat{\mathbf{u}}(z_k)$ with $\phi_{k}=kh$.
If the contour $\Gamma$ is symmetric with respect to the real axis, $A$ is a real matrix, then $\widehat{u}(\overline{z})=\overline{\widehat{u}(z)}$, and after truncation, there is
\begin{equation}\label{eq:finapprox}
  \mathbf{u}(t)\approx \frac{h}{\pi} \textmd{Im} \left\{\sum\limits_{k = 0}^{N-1}e^{z_{k}t}\mathbf{\widehat{u}}_{k}z_{k}'\right\}.
\end{equation}
This is the CIM scheme of (\ref{eq:iniproblem1}).

The key to design an efficient CIM scheme is to find the spectrum distribution of the matrix $A$, which determines how to choose an appropriate contour integral curve $\Gamma$.

The non-locality and the weakly singular kernel of the time fractional operator results to the $\mathscr{O}\left(N^2\right)$ computational cost for time-marching scheme and weak singularity of the solution. A lot of efforts have been made to efficiently deal with this difficulty (see, e.g., \cite{SIAM18,Freed02} etc).
Compared with the time-stepping methods (see eg. \cite{Freed02,Zeng22}), the CIM scheme has the following cons and pros, when solving the nonlocal problems.
\begin{itemize}
  \item Generally, the time-stepping methods need the memory of $\mathscr{O}(N)$ and have the computational complexity of  $\mathscr{O}\left(N^2\right)$. While for the CIM scheme, the required memory is $\mathscr{O}(1)$ and  the computational complexity is $\mathscr{O}(N)$.

  \item  For the time-stepping methods, the solution at a given  later time depends on the previous ones. While, for the CIM scheme, the solution can be directly computed at any desired time, without the information on earlier time.

  \item The computation cost of the CIM scheme mainly lies in solving the system (\ref{eq:Lapalce}), which can be parallelly computed with the rate of $100\%$.

  \item For the time-stepping methods, low regularity of the solution will make it hard to get a high convergence rate.  This issue has little influence on the CIM scheme.

  \item[$(\bullet)$] There is nothing perfect. Although, the CIM works well for the linear model, it is difficult to deal with the nonlinear one directly.

\end{itemize}

Above all, the CIM is a simple, time-saving, and efficient numerical method. To build a CIM scheme, the key is to choose an appropriate integral contour, which depends on the spectrum distribution of the matrix $A$. Currently, there are four types of popular integral contours used for the CIM, namely, Talbot's contour \cite{Talbot79}, parabolic contour e.g., \cite{Weidemantrefe07},  hyperbolic contour  e.g.,  \cite{Fernandez06,Ma22}, and other simple, closed, and positively oriented curves, e.g., \cite{Sloan99,Dingfelder15}. The CIM with these contours can be used to solve parabolic problems, e.g., \rred{\cite{Sloan99,Lee06,Li21}}, integral differential equation with convolution memory kernel \cite{Sheen06,Li22,Li23}, Black-Scholes and Heston equations \cite{int11}, and other problems, e.g., \cite{XiYang18,ZhouSun18}. During these applications, the CIM behaves high numerical performance. This paper will develop the CIM scheme into the time fractional differential system, i.e., the Feynman-Kac equation with two internal states \cite{xu18}.

Feynman-Kac equation usually describes the distribution of a particular type of statistical observables, e.g., functional of the particle trajectory \cite{SLE_07,wudeng16,xu18}.  The model considered in this paper characterizes a specific functional: $A=\int_{0}^{t}U(j(\tau))d \tau$, where $j(\tau)$ represents the $j$-th internal state at time $\tau$ with values belonging to $\{1, 2, \cdots, m\}$. The distribution of $A$ in the frequency domain is governed by
\begin{equation}\label{eq:Feynmankac}
\left\{
\begin{aligned}
\mathbf{M}^{T}\frac{\partial}{\partial t}\mathbf{G}&=\left(\mathbf{M}^{T}-\mathbf{I}\right)
       \mathrm{diag}\left(B_{\alpha_{1}}^{-1}, B_{\alpha_{2}}^{-1},\cdot\cdot\cdot,B_{\alpha_{m}}^{-1}\right)
       \mathrm{diag}\left(\mathfrak{D}_{t}^{1-\alpha_{1}}, \mathfrak{D}_{t}^{1-\alpha_{2}}, \cdot\cdot\cdot,
       \mathfrak{D}_{t}^{1-\alpha_{m}}\right)\mathbf{G}\\
       &\quad -\rho\mathbf{M}^{T} \mathrm{diag}\left(U(1),U(2),\cdot\cdot\cdot, U(m)\right)\mathbf{G},\quad t\in(0,T],\\
\mathbf{G}(\cdot,0) &= \mathbf{G_{0}},
\end{aligned}\right.
\end{equation}
where $M$ is the transition matrix of a Markov chain with dimension $m\times m$; $B_{\alpha_{j}}^{-1}$, $j=1,2,\cdot\cdot\cdot,m$, are given positive real numbers;
$\mathbf{G}=[G_{1},G_{2},\cdot\cdot\cdot,G_{m}]^{T}$ denotes the solution of the model (\ref{eq:Feynmankac}) with $G_{j}:= G_{j}(\rho,t)$ represents the Laplace transform of $G_{j}(A,t)$ w.r.t. $A$; and $G_{j}(A,t)$ is the probability density function (PDF) of finding the particle with the functional $A$ in the $j$-th internal state at time $t$;
$I$ is the identity matrix; `$\mathrm{diag}$' represents the diagonal matrix formed by its vector arguments; and $\mathfrak{D}_{t}^{1-\alpha_{j}}$, $j=1,2,\cdot\cdot\cdot,m$, are the fractional substantial derivatives, defined as
\begin{small}\begin{equation} \label{DefinationFSD}
 \mathfrak{D}^{1-\alpha_{j}}_{t}G_{j}(\rho,t):=\left(\frac{\partial}{\partial t}
 +\rho U(j)\right)\frac{1}{\Gamma(\alpha_{j})}\int_{0}^{t}\frac{\mathrm{exp}\left[-(t-\tau)\rho U(j)\right]}{(t-\tau)^{1-\alpha_{j}}}G_{j}(\rho,\ \tau)d\tau
\end{equation}\end{small}
with $0<\alpha_{j}<1, \,j=1,2,\cdot\cdot\cdot,m$.

This paper is organized as follows. In Section \ref{sec:Preliminaries}, we give the regularity estimates on the solution of (\ref{eq:Feynmankac}). In Section \ref{subsec:NCIM}, the CIMs with parabolic contour and hyperbolic contour for the system (\ref{eq:Feynmankac}) are built respectively. Also we perform the error estimate and stability analysis for these schemes. In addition, the parameters in parabolic and hyperbolic contours are optimally determined. To verify the efficiency of the CIMs, we also construct a time-marching scheme To provide a reference solution. Some numerical experiments are performed in Section \ref{sec:NumericalResults} to show the high numerical performance of the CIMs in solving such a non-local system. Concluding remarks are presented in Section \ref{sec:conclusions}.

\section{The continuous problem}
\label{sec:Preliminaries}
In this section, we perform the regularity analysis for Problem (\ref{eq:Feynmankac}).
\subsection{Solution representations}
\label{subsec:twointernalstatesFeynman-Kac}
We consider the Feynman-Kac equation with two internal states. Without loss of generality, the transition matrix $M$ can be written as
\begin{equation}\label{eq:matrices}
 M=\begin{bmatrix}p&1-p\\1-b&b\end{bmatrix},
\end{equation}
where $p$, $b\in[0, 1]$, and $p+b\neq1$. Then Problem (\ref{eq:Feynmankac}) reduces to
\begin{equation}\label{eq:DIVFeynman-kac}
 \left\{
 \begin{aligned}
&p\left(\frac{\partial}{\partial t}+U(1)\rho\right)G_{1}+(1-b)\left(\frac{\partial}{\partial t}+U(2)\rho\right)G_{2}
=(p-1)B_{\alpha_{1}}^{-1}\mathfrak{D}_{t}^{1-\alpha_{1}}G_{1}+(1-b)B_{\alpha_{2}}^{-1}\mathfrak{D}_{t}^{1-\alpha_{2}}G_{2},\ t\in(0,T],\\
&(1-p)\left(\frac{\partial}{\partial t} + U(1)\rho\right)G_{1} + b\left(\frac{\partial}{\partial t}+U(2)\rho\right)G_{2}
=(1-p)B_{\alpha_{1}}^{-1}\mathfrak{D}_{t}^{1-\alpha_{1}}G_{1}+(b-1)B_{\alpha_{2}}^{-1}\mathfrak{D}_{t}^{1-\alpha_{2}}G_{2},\ t\in(0,T],\\
&\mathbf{G}(\cdot,0)= \mathbf{G_{0}},\ \quad\quad\quad\quad\quad\quad\quad\quad\quad\quad\quad&
\end{aligned}\right.
\end{equation}
where $\mathbf{G}_{0}=[G_{1,0},G_{2,0}]^{T}$ is the initial value.
After simple calculations, we can obtain from (\ref{eq:DIVFeynman-kac}) that
\begin{displaymath}
\left(\frac{\partial}{\partial t} + U(1)\rho\right)G_{1}+\left(\frac{\partial}{\partial t} + U(2)\rho\right)G_{2}=0.
\end{displaymath}
Then, there is
\begin{equation}\label{bianxing}
 \left\{
 \begin{aligned}
\left(\frac{\partial}{\partial t} + U(1)\rho\right)G_{1}
&=\frac{1-p}{1-p-b} B_{\alpha_{1}}^{-1}\mathfrak{D}_{t}^{1-\alpha_{1}}G_{1}+\frac{1-b}{p+b-1}
  B_{\alpha_{2}}^{-1}\mathfrak{D}_{t}^{1-\alpha_{2}}G_{2},\\
\left(\frac{\partial}{\partial t} + U(2)\rho\right)G_{2}
&=\frac{1-p}{p+b-1} B_{\alpha_{1}}^{-1}\mathfrak{D}_{t}^{1-\alpha_{1}}G_{1}+\frac{1-b}{1-p-b}
  B_{\alpha_{2}}^{-1}\mathfrak{D}_{t}^{1-\alpha_{2}}G_{2}.
\end{aligned}\right.
\end{equation}
Denote $m_{1}:=\frac{1-p}{1-p-b}$, $m_{2}:=\frac{1-b}{1-p-b}$. For $0<\alpha_{1}, \alpha_{2}<1$, if $G_j\in I_{0^+}^{1-\alpha_j}[L_{1}(0,T)]:=\{f:f(x)={}_{0}I^{1-\alpha_j}_{t}\varphi(t), \varphi(t)\in L_{1}(0,T)\}$, then ${}_{0}D_{t}^{-\alpha_{j}}(e^{\rho U(j)t}G_{j})\big|_{t=0}=0$, $j=1, 2$ (we note that the space $I_{0^+}^{1-\alpha_j}[L_{1}(0,T)]$ only excludes some extreme functions such as
 $t^{-\beta}$, $\beta<\alpha_j$, so it is not a harsh requirement for $G_j$, see more details in \cite{zhao21}).  By this, taking the Laplace transform on both sides of (\ref{bianxing}), we deduce
\begin{equation}\label{eq:twointernal}
 \left\{\begin{aligned}
\widehat{G}_{1}&=\left(\left(z+\rho U(1)\right)-m_{1}B_{\alpha_{1}}^{-1}\left(z+\rho U(1)\right)^{1-\alpha_{1}}\right)^{-1}
               \left(-m_{2}B_{\alpha_{2}}^{-1}\left(z+\rho U(2)\right)^{1-\alpha_{2}}\widehat{G}_{2}+G_{1,0}\right),\\
\widehat{G}_{2}&=\left(\left(z+\rho U(2)\right)-m_{2}B_{\alpha_{2}}^{-1}\left(z+\rho U(2)\right)^{1-\alpha_{2}}\right)^{-1}
               \left(-m_{1}B_{\alpha_{1}}^{-1}\left(z+\rho U(1)\right)^{1-\alpha_{1}}\widehat{G}_{1}+G_{2,0}\right).
\end{aligned}\right.
\end{equation}
Let
\begin{equation}\label{eq:Hfunction}
H(z):=\left\{[(z+\rho U(1))^{\alpha_{1}}-m_{1}B_{\alpha_{1}}^{-1}][(z+\rho U(2))^{\alpha_{2}}-m_{2}B_{\alpha_{2}}^{-1}]- m_{1}m_{2}B_{\alpha_{1}}^{-1}B_{\alpha_{2}}^{-1}\right\}^{-1},
\end{equation}
\begin{equation}\label{eq:Halpha1}
H_{\alpha_{1}}(z):=H(z)\left((z+\rho U(2))^{\alpha_{2}}-m_{2}B_{\alpha_{2}}^{-1}\right),
\end{equation}
and
\begin{equation}\label{eq:Halpha2}
H_{\alpha_{2}}(z):=H(z)\left((z+\rho U(1))^{\alpha_{1}}-m_{1}B_{\alpha_{1}}^{-1}\right).
\end{equation}
Then, the system (\ref{eq:DIVFeynman-kac}) can be decoupled in Laplace space as
\begin{equation}\label{eq:solutionsinlaplace}
 \left\{\begin{aligned}
\widehat{G}_{1}&=
               H_{\alpha_{1}}(z)\big(z+\rho U(1)\big)^{\alpha_{1}-1}G_{1,0}-m_{2}B_{\alpha_{2}}^{-1}H(z)\big(z+\rho U(1)\big)^{\alpha_{1}-1}G_{2,0},\\
\widehat{G}_{2}&=
               H_{\alpha_{2}}(z)\big(z+\rho U(2)\big)^{\alpha_{2}-1}G_{2,0}-m_{1}B_{\alpha_{1}}^{-1}H(z)\big(z+\rho U(2)\big)^{\alpha_{2}-1}G_{1,0}.
\end{aligned}\right.
\end{equation}

\subsection{Regularity}\label{subsec:Regularity}
Define the sectors
\begin{displaymath}
\Sigma_{\theta}:=\left\{z\in\mathbb{C}: |\arg(z)|\leq \theta, z\neq0\right\},\quad\theta\in(\pi/2,\pi);
\end{displaymath}
\begin{displaymath}
\Sigma_{\theta,\delta}:=\left\{z\in\mathbb{C}: |z|>\delta>0,|\arg(z)|\leq \theta\right\},\quad\theta\in(\pi/2,\pi).
\end{displaymath}
Take an integral contour 
\begin{displaymath}
\Gamma_{\theta,\delta}:=\{z\in\mathbb{C}:|z|=\delta>0,|\arg(z)|\leq \theta\}\cup\left\{z\in\mathbb{C}:z=re^{\pm i\theta},r\geq\delta>0\right\},
\end{displaymath}
oriented with an increasing imaginary part. Based on these analytic settings, we have the following estimates related to (\ref{eq:DIVFeynman-kac}) or (\ref{eq:solutionsinlaplace}). See \ref{subsec:RegularityProof} for the proofs in details.
\begin{lem}\label{lemm:lemma2.4}
Let $z\in\Sigma_{\theta,\delta}$ and $\delta\geq 2 \max\limits_{j=1,2}\left\{|\rho U(j)|\right\}$. For $0<\alpha_{j}<1$, there hold
\begin{displaymath}
\left|(z+\rho U(j))^{-\alpha_{j}}\right|\leq 2 |z|^{-\alpha_{j}},\ \  j=1,2.
\end{displaymath}
\end{lem}

\begin{lem}\label{lemm:lemma2.5}
Let $z\in\Sigma_{\theta,\delta}$ and $\delta\geq \max\limits_{j=1,2}\left\{2 |\rho U(j)|,  4^{1/\alpha_j}|m_{j}B_{\alpha_{j}}^{-1}|^{1/\alpha_{j}}\right\}$. For $0<\alpha_{j}<1$, there hold
\begin{displaymath}
\left|((z+\rho U(j))^{\alpha_{j}}-m_{j}B_{\alpha_{j}}^{-1})^{-1}\right| \leq 4 |z|^{-\alpha_{j}},\ \  j=1,2.
\end{displaymath}
\end{lem}

\begin{lem}\label{lemm:lemma2.6}
Let $z\in\Sigma_{\theta,\delta}$, $\delta\geq\max\left\{
\max\limits_{j=1,2}\left\{2 |\rho U(j)|, 4^{1/\alpha_j}|m_{j}B_{\alpha_{j}}^{-1}|^{1/\alpha_{j}}\right\},  \left(32|m_{1}m_{2}B^{-1}_{\alpha_{1}}B^{-1}_{\alpha_{2}}|\right)^{1/(\alpha_{1}+\alpha_{2})}\right\}$. For $0<\alpha_{j}<1$, there hold
\begin{displaymath}
\left|H(z)\right|\leq 32 |z|^{-\alpha_{1}-\alpha_{2}},\
\left|H_{\alpha_{j}}(z)\right|\leq  8 |z|^{-\alpha_{j}},~~j=1,2,
\end{displaymath}
where $H(z)$ and $H_{\alpha_{j}}(z)$ are defined in (\ref{eq:Hfunction}), (\ref{eq:Halpha1}), and (\ref{eq:Halpha2}),
 respectively. 
\end{lem}

\begin{lem}\label{lemm:lemma2.8}
Let $z\in\Sigma_{\theta,\delta}$ and $\delta$ satisfy the conditions in Lemma \ref{lemm:lemma2.6}. For $0<\alpha_{j}<1$, there hold
\begin{displaymath}
\left|H_{\alpha_{j}}(z)(z+\rho U(j))^{\alpha_{j}-1}\right|\leq 16 |z|^{-1},
\end{displaymath}
and
\begin{displaymath}
\left|H(z)(z+\rho U(j))^{\alpha_{j}-1}\right|\leq 64 |z|^{-1-\alpha_{j^{*}}},~~j=1,2,~~j\cdot j^{*}=2.
\end{displaymath}
\end{lem}

Based on the above results, one can get the estimates on the solutions of (\ref{eq:DIVFeynman-kac}).
\begin{thm}\label{the:theorem2.11}
Under the conditions in Lemma \ref{lemm:lemma2.6}, for given initial values $G_{j,0}$, $j=1,2,$ and $t>1/\delta$, the solutions of (\ref{eq:DIVFeynman-kac}) satisfies
\begin{displaymath}
\left|G^{(q)}_{j}(t)\right|\leq
      \frac{32\left|m_{j^{*}}B_{\alpha_{j^{*}}}^{-1}\right|}{\pi}\left((-1/\cos(\theta))t^{\alpha_{j^{*}}-q}e^{\cos(\theta)}+2 \theta\delta^{q-\alpha_{j^{*}}}e^{\delta T}\right)|G_{j^{*},0}|
      +\frac{8}{\pi}\left((-1/\cos(\theta))t^{-q}e^{\cos(\theta)}+2\theta\delta^{q}e^{\delta T}\right)|G_{j,0}|,
\end{displaymath}
for $q=0,1$, where $G^{(q)}_{j}(t)$ denote $\frac{\partial^{q}}{\partial t^{q}}G_{j}(t)$, $j=1,2,$ and $j\cdot j^{*}=2$.
\end{thm}

Theorem \ref{the:theorem2.11} discovers the weak singularity of the solution bear to the origin, which usually weakens the performance of the time-marching schemes but, we will see in the following context, has no influence on the CIMs.

\section{The schemes and error estimates for the Feynman-Kac system}
\label{subsec:NCIM}
In this section, for the system (\ref{eq:DIVFeynman-kac}), the CIMs with two kinds of contours, i.e., parabolic contour and hyperbolic contour, are given. With careful analysis of the analytical domain of the solution in frequency domain, the parameters used in the contours are determined. The error estimates and stability analysis are presented. To give the reference solution for verifying the effectiveness of the CIMs, a time marching scheme is also designed.
\subsection{The CIM schemes}
Here we discuss two different integral contours \cite{Weidemantrefe07} for the CIMs, which are parameterized by
\begin{equation}\label{eq:paraboliccontour}
 \mathrm{( Parabolic \ contour)}\ \ \ \Gamma_{1}:\ z(\phi) = \eta_{1}(i\phi+1)^{2}, ~ -\infty<\phi<\infty,
\end{equation}
and
\begin{equation}\label{eq:hyperboliccontour}
\mathrm{( Hyperbolic \ contour)}\ \ \ \Gamma_{2}:\ z(\phi) = \eta_{2}(1+\sin(i\phi-\alpha)), ~ -\infty<\phi<\infty,
\end{equation}
where $\eta_{1}$, $\eta_{2}>0$ and $\alpha>0$ are the parameters to be determined. With these, the solutions $G_{1}(t)$, $G_{2}(t)$ can be represented as the infinite integrals with respect to $\phi$, i.e.,
\begin{equation}\label{eq:integralsolve}
G_{j}(t)=I_{j}:=\int_{-\infty}^{+\infty} v_{j}(t,\phi)d\phi,\quad j=1, 2,
\end{equation}
where
\begin{equation}\begin{aligned}\label{exp:V1V2}
v_{j}(t,\phi):=\frac{1}{2\pi i}e^{z(\phi)t}\widehat{G}_{j}(z(\phi))z'(\phi), ~~ j=1,2 
\end{aligned}\end{equation}
with $\widehat{G}_{j}$ defined in (\ref{eq:twointernal}).

Applying the trapezoidal rule to compute the integral (\ref{eq:integralsolve}) with uniform steps $h_{m}$ $(m=1, 2)$ for $\phi$, one can numerically get the approximate solutions $G^{\langle m\rangle}_{1,N}(t)$ and $G^{\langle m\rangle}_{2,N}(t)$. 
Since the contours $\Gamma_{1}$ and $\Gamma_{2}$ are symmetric with respect to the real axis and $\widehat{G}_{j}(\overline{z(\phi)})=\overline{\widehat{G}_{j}(z(\phi))}$, $j=1,2$, (\ref{eq:integralsolve}) can be approximated by
 \begin{equation}\begin{aligned}\label{eq:computG1}
  G^{\langle m\rangle}_{j,N}(t)\approx I_{j,h_{m},N}:= h_{m}\sum\limits_{k = 1-N}^{N-1}v_{j}(t,\phi_{k})=\frac{h_{ m}}{\pi}\mathrm{Im}\left\{\sum\limits_{k = 0}^{N-1}e^{z_{k}t}\widehat{G}_{j,k}z'_{k}\right\},\ \ j=1,2,
 \end{aligned}\end{equation}
where $m=1,2$ denote the two different choices of the integral contours $\Gamma_{1}$ and $\Gamma_{2}$.

It can be seen from (\ref{eq:computG1}) that the numerical solutions at current time are computed without knowing its information at previous times.

\subsection{Quadrature error}
\label{sec:err}

The key issue to ensure the effectiveness and efficiency of the CIMs is to make all singular points of $\widehat{G}_{j}$ locate at the left side of $\Gamma_{m}$ and that the integrand of the indefinite integrals (\ref{eq:integralsolve}) has wide analytical open strips see, e.g., \cite{Weidemantrefe07,Trefethren}.

\subsubsection{Determination of the open strip}
\label{sec:openstrip}

Since $z(\phi)$ in both of the integral contours (\ref{eq:paraboliccontour}) and (\ref{eq:hyperboliccontour}) are analytical,  the analytical properties of $v_{j}(t,\phi)$ in (\ref{exp:V1V2}) is completely determined by $\widehat{G}_{j}(z(\phi))$. According to the expression of $\widehat{G}_{j}(z(\phi))$ in (\ref{eq:solutionsinlaplace}),  the singularity mainly comes from $H(z)$.

Here we consider the case $U(j)\geq0$, and denote $C_{j}:=m_{j}B_{\alpha_{j}}^{-1}$, $j=1,2$. There are $C_{j}\neq0$, $j=1,2$. The analytical domain of $H(z)$ are determined by the next proposition.

\begin{prop}\label{Thm:Prop1}
Let $z\in\sum_{\theta}$. If
\begin{equation}\label{condition01}
\textmd{Re}\left((z+\rho U(j))^{\alpha_j}\right)>2 C_j,~~j=1,2,
 ~~ {\rm or } ~~~
\left|\textmd{Im}\left((z+\rho U(j))^{\alpha_j}\right)\right|>|C_j|,~~j=1,2,
\end{equation}
then $H(z)$ is analytic.
\end{prop}

See the proof in  \ref{det:strip}.

\begin{rmk}\label{remark_d1_d2}
From (\ref{condition01}), it can be further obtained that $H(z)$ is analytic if
\begin{equation}\label{condition02}
{\rm Re}(z)>d_1 ~~ {\rm or }~~ \left|{\rm Im}(z)\right|>d_2,
\end{equation}
where $d_1:=\max\limits_{j=1,2}\left\{(2|C_j|/\cos(\alpha_j\pi/2))^{1/\alpha_j}-\textmd{Re}(\rho U(j))\right\}$ and $d_2:=\max\limits_{j=1,2}\left\{(|C_j|/\cos(\alpha_j\pi/2))^{1/\alpha_j}-\textmd{Im}(\rho U(j))\right\}$. One can see \ref{condi} for more details.
\end{rmk}

\begin{figure}[h]
\centering
\begin{minipage}[c]{0.40\textwidth}
 \centering
 \centerline{\includegraphics[height=4.2cm,width=5.7cm]{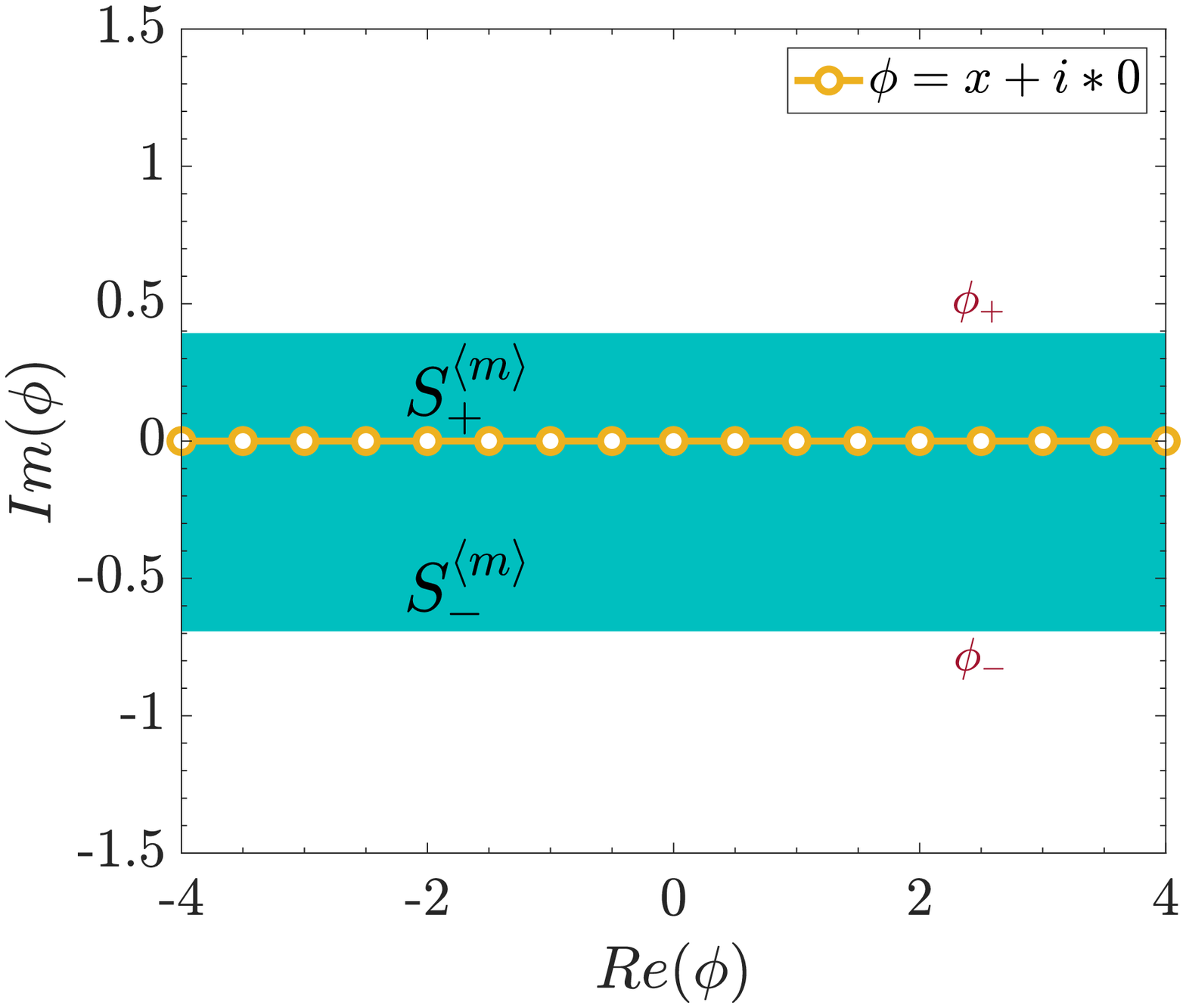}}
\end{minipage}
\begin{minipage}[c]{0.40\textwidth}
 \centering
 \centerline{\includegraphics[height=4.2cm,width=5.7cm]{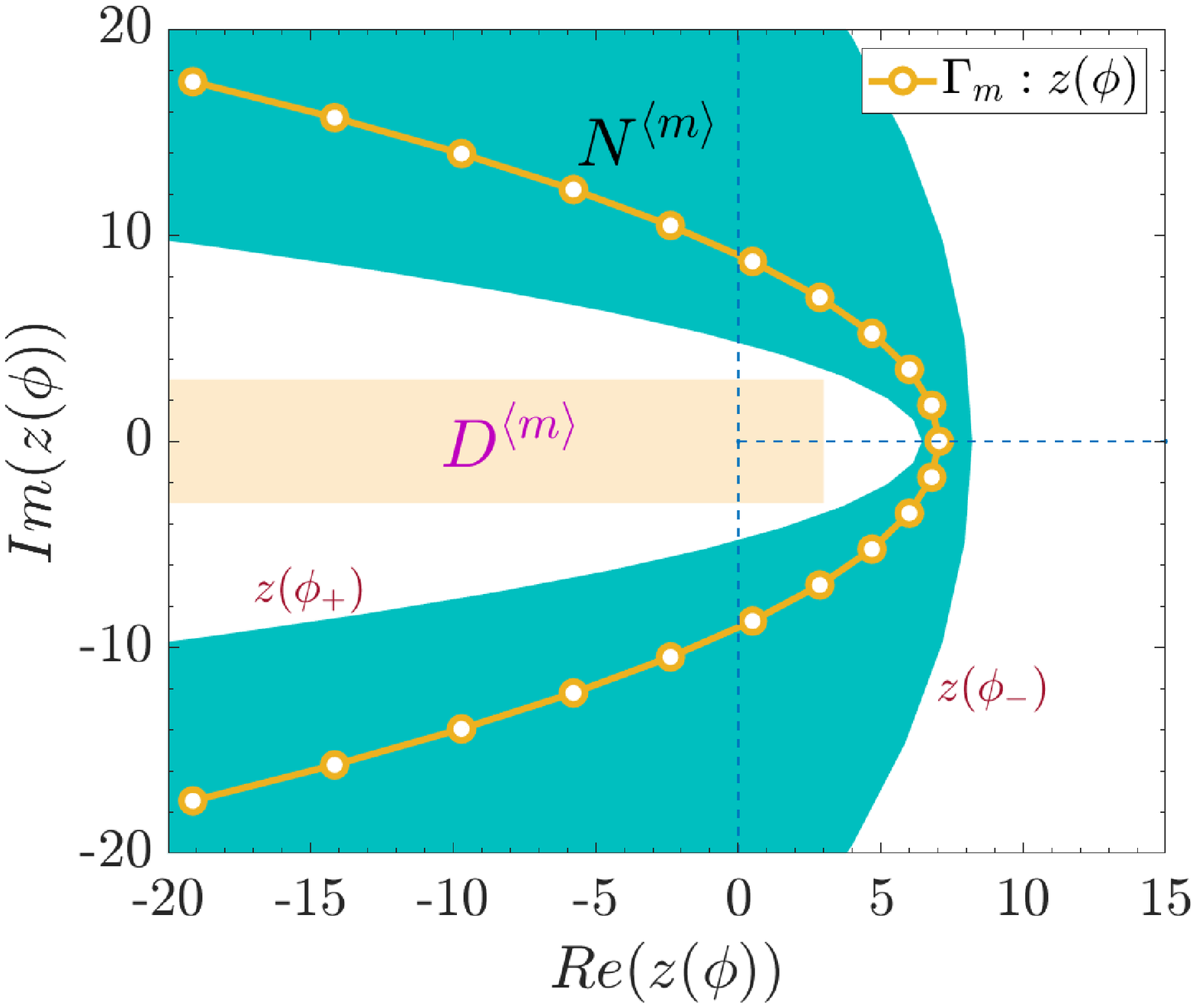}}
\end{minipage}
\caption{The schematic diagrams of open strip $S^{\langle m\rangle}$ and the corresponding neighbourhood $N^{\langle m\rangle}$; (\emph{left}) is the horizontal open strip $S^{\langle m\rangle}$; (\emph{right}) is the neighborhood $N^{\langle m\rangle}$ obtained by mapping the horizontal strip $S^{\langle m\rangle}$ by the conformal transformation $z(\phi)$.}
\label{Fig:strip}
\end{figure}

Given open strip
\begin{displaymath}
S^{\langle1\rangle}:=\{\phi=x+iy\in\mathbb{C}: -d<y<a<1, x\in\mathbb{R}~{\rm and}~a, d>0\}
\end{displaymath}
or
\begin{displaymath}
S^{\langle2\rangle}:=\{\phi=x+iy\in\mathbb{C}: -\alpha<y<\pi/2-\alpha-\delta, x\in\mathbb{R}, 0<\alpha, \delta<\pi/2~{\rm and}~\alpha<\pi/2-\delta\},
\end{displaymath}
then the integral contour $\Gamma_1$ or $\Gamma_2$ maps it into a neighbourhood, defined as $N^{\langle m\rangle}:=\{z(\phi):\phi\in S^{\langle1\rangle}~ {\rm or}~ S^{\langle2\rangle}\}$. According to (\ref{condition02}), there exists a domain $D^{\langle m\rangle}$, such that $H(z)$ is analytical in the domain $\mathbb{C}\setminus D^{\langle m\rangle}$, shown in Figure \ref{Fig:strip} (\emph{right}). Hence, for any $z(\phi)\in N^{\langle m\rangle}\subset\mathbb{C}\setminus D^{\langle m\rangle}$, the integrand $v_j(t, \phi)$ is analytical. Since the integral contours $\Gamma_1$ and $\Gamma_2$ are holomorphic mappings, so, once $N^{\langle m\rangle}$ are specified, the strips $S^{\langle1\rangle}$ and $S^{\langle2\rangle}$ can be accordingly determined.


For the consideration of convergence, we need to select the analytical domain $N^{\langle m\rangle}$ as large as possible. Correspondingly, the optimum parameter $a$ in $S^{\langle1\rangle}$ is $a=1-\left[(d_{1}+(d_{1}^{2}+d_{2}^{2})^{1/2})/(2\eta_{1})\right]^{1/2}$ and the biggest value of $d$ is given by (\ref{eq:d}); the optimal parameter $\delta$ in $S^{\langle2\rangle}$ is $\delta=\arctan\left(d_{2}/(\eta_{2}-d_{1})\right)$ and $\alpha$ is determined by maximizing $Q(\alpha)=\frac{\pi^{2}-2\pi\alpha-2\pi\delta}{A(\alpha)}$ (see Subsection \ref{sec:opt}). For more details, one can see \ref{det:strips}.

\subsubsection{Stability}

Paved by the front part, the integrands $v_{j}(t, \phi)$ are analytical with respect to $\phi$ in the open strips $S^{\langle 1\rangle}$ and $S^{\langle 2\rangle}$. Now, for the CIMs (\ref{eq:computG1}) with the integral contours $\Gamma_{1}$ and $\Gamma_{2}$, we have the following stability analyses.

\begin{lem}[\cite{Fernandez04}, Lemma 1]\label{et:lemma}
Take $L(x):=1+\big|\ln(1-e^{-x})\big|$, $x>0$, there hold
\begin{displaymath}
\int_{0}^{+\infty}e^{-\gamma\cosh(x)}dx\leq L(\gamma),\ \gamma>0,
\end{displaymath}
and
\begin{displaymath}
\int_{\sigma}^{+\infty}e^{-\gamma\cosh(x)}dx\leq (1+L(\gamma))e^{-\gamma\cosh(\sigma)},\ \gamma>0\ {\rm and}~ \sigma>0.
\end{displaymath}
\end{lem}

\begin{lem}\label{lemm:lemma3.1}
Let $v_{1}(t,\phi)$ be defined in (\ref{exp:V1V2}) and $z(\phi)$ defined in (\ref{eq:paraboliccontour}) with $\phi=x+ir \in S^{\langle 1\rangle}$.
For $t>0$ and $B>0$ defined in (\ref{con:B}), we have
\begin{equation}\label{est:bound1}
\left|v_{1}(t,\phi)\right| \leq \frac{B e^{\eta_{1}t}}{\pi}\left(|G_{1,0}|
+\eta_{1}^{-\alpha_{2}}|G_{2,0}|\right)e^{-\eta_{1}tx^{2}} ~~ \forall~\phi\in S^{\langle 1\rangle}.
\end{equation}
Let $z(\phi)$ defined in (\ref{eq:hyperboliccontour}) with $\phi=x+iy \in S^{\langle 2\rangle}$. For $t>0$ and $0<y<\min\{\alpha, \pi/2-\alpha-\delta\}$, there holds
\begin{equation}\label{est:bound2}
\left|v_{1}(t,\phi)\right|
\leq \frac{8e^{\eta_{2}t}}{\pi}\sqrt{\frac{1+\sin (\alpha+y)}{1-\sin (\alpha+y)}}\left(|G_{1,0}|
     +\frac{4\left|m_1B_{\alpha_1}^{-1}\right|}{\eta_{2}^{\alpha_{2}}(1-\sin (\alpha+y))^{\alpha_{2}}}|G_{2,0}|\right)
     e^{-\eta_{2}t\sin{\alpha}\cosh x}~ ~ \forall\phi\in S^{\rred{\langle 2\rangle}}.
\end{equation}

For the integrand $v_{2}(t,\phi)$, similar estimates can be obtained.
\end{lem}
The proof of  Lemma \ref{lemm:lemma3.1} is given in \ref{proof:lemm3}.

It can be seen from Lemma \ref{lemm:lemma3.1} that the decay of the solution mainly depends on the size of the real part.
The stability results of the CIMs are given as follows.

\begin{thm}\label{thm:estev}
Let $G^{\langle 1\rangle}_{1,N}(t)$ be defined in (\ref{eq:computG1}) with uniform steps $h_1$. For $t>0$ and $B>0$ defined in (\ref{con:B}), there holds
\begin{displaymath}
\left|G^{\langle 1\rangle}_{1,N}(t)\right|\leq B\sqrt{\frac{\eta_{1}}{\pi}}\left(\eta_{1}^{-1}|G_{1,0}|
+\eta_{1}^{-(1+\alpha_{2})}|G_{2,0}|\right)t^{-1/2}e^{\eta_{1}t}.
\end{displaymath}
Let $G^{\langle 2\rangle}_{1,N}(t)$  be defined in (\ref{eq:computG1}) with uniform steps $h_2$. For $t>0$ and $0<y<\min\{\alpha, \pi/2-\alpha-\delta\}$, there holds
\begin{displaymath}
\left|G^{\langle 2\rangle}_{1,N}(t)\right|
\leq \frac{16}{\pi}\sqrt{\frac{1+\sin
    (\alpha+y)}{1-\sin (\alpha+y)}}\left(|G_{1,0}|
+\frac{4\left|m_1B_{\alpha_1}^{-1}\right|}{\eta_{2}^{\alpha_{2}}(1-\sin
    (\alpha+y))^{\alpha_{2}}}|G_{2,0}|\right)L(\eta_{2}t\sin(\alpha))e^{\eta_{2}t}.
\end{displaymath}

For $\left|G^{\langle m\rangle}_{2,N}(t)\right|$, $m=1,2$, similar estimates can be obtained.
\end{thm}

The proof of Theorem \ref{thm:estev} can be found in \ref{proof:thm}. Notice that $L(x)$ given in Lemma \ref{et:lemma} is decreasing; $L(x)\rightarrow1$  as $x\rightarrow\infty$ and $L(x)\sim|\ln x|$ as $x\rightarrow0^{+}$. It can be seen that the CIMs are unconditionally stable with respect to the initial values.

\subsubsection{Error estimates and determination of the optimal parameters}
\label{sec:opt}
Here, we will determine the optimal parameters in the parabolic integral contour $\Gamma_{1}$ and hyperbolic contour $\Gamma_{2}$, respectively, and prove that  the CIMs (\ref{eq:computG1}) of the Feynmann-Kac equation with two internal states have spectral accuracy.

Denote $I_{j,h_m}:=h_{m}\sum_{k=-\infty}^{\infty}v_{j}(t,kh_m)$. Then the error of the CIMs can be expressed as
\begin{displaymath}
E^{\langle m\rangle}_{j,N}:=\left|I_{j}-I_{j,h_m,N}\right|\leq DE^{\langle m\rangle}+TE^{\langle m\rangle},~ ~ m,j=1,2,
\end{displaymath}
where $TE^{\langle m\rangle}=\left|I_{j,h_m}-I_{j,h_m,N}\right|$ is the truncation error and $DE^{\langle m\rangle}:=DE^{\langle m\rangle}_{+}+DE^{\langle m\rangle}_{-}:=|I_{j}-I_{j,h_m}|$ is the discretization error. The standard estimates of the  discretization error are shown in the following lemma.

\begin{lem}[\cite{Weidemantrefe07} Theorem 2.1]\label{lem:lemma2}
Consider the absolutely convergent integral
$$I:=\int_{-\infty}^{\infty}g(u)du,$$ and its infinite and trapezoidal approximation
$$I_\tau:=\tau\sum_{k=-\infty}^{\infty}g(k\tau).$$ Let $w=u+iv$, with $u$ and $v$ real. Suppose $g(w)$ is analytic in the strip $-d<v<c$, for some $c>0$, $d>0$, with $g(w)\rightarrow0$ uniformly as $|w|\rightarrow\infty$ in this strip. Suppose further that for some $M_+>0$, $M_->0$, the function $g(w)$ satisfies
\begin{equation}\label{eq:MMM}
\int_{-\infty}^{\infty}|g(u+ir)|du\leq M_+,~ \int_{-\infty}^{\infty}|g(u-is)|du\leq M_-,
\end{equation}
for all $0<r<c$, $0<s<d$.
Then the discretization error
\begin{displaymath}
DE:=|I-I_\tau|\leq DE_{+}+DE_{-},
\end{displaymath}
where
\begin{displaymath}
DE_{+}=\frac{M_{+}}{e^{2\pi c/\tau}-1}, ~ ~ DE_{-}=\frac{M_{-}}{e^{2\pi d/\tau}-1}.
\end{displaymath}
\end{lem}

Based on lemma \ref{lem:lemma2}, for the CIMs with the integral contour $\Gamma_{1}$ and $\Gamma_{2}$, we have the following error estimates.

\begin{thm}\label{thm:accuracy}
Let $G_1(t)$ and $G^{\langle 1\rangle}_{1,N}(t)$ be the solutions of (\ref{eq:DIVFeynman-kac}) and (\ref{eq:computG1}) with uniform step-size $h_1$. Given $a>0$, for $t>0$, $h_1=\mathscr{O}(1/N)$, there holds
\begin{displaymath}
\left|G_1(t)-G^{\langle 1\rangle}_{1,N}(t)\right|
\leq B\sqrt{\frac{\eta_{1}}{\pi}}t^{-1/2}\left(\eta_{1}^{-1}|G_{1,0}|
  +\eta_{1}^{-(1+\alpha_{2})}|G_{2,0}|\right)
\left(\frac{e^{\eta_{1}t}}{e^{2\pi a/h_{1}}-1}+\frac{e^{\pi^{2}/(\eta_{1}th_{1}^{2})}} {e^{2\pi^{2}/(\eta_{1}th_{1}^{2})-2\pi/h_1}-1}+\frac{e^{\eta_{1}t}}{e^{\eta_{1}t(Nh_{1})^{2}}}\right).
\end{displaymath}
Let $G^{\langle 2\rangle}_{1,N}(t)$ be the solution of (\ref{eq:computG1})  with uniform step-size $h_2$. For $t>0$, $h_2=\mathscr{O}(1/N)$, and $0<y<\min\{\alpha, \pi/2-\alpha-\delta\}$, we have
\begin{displaymath}\begin{aligned}
\quad \left|G_1(t)-G^{\langle 2\rangle}_{1,N}(t)\right|
&\leq \frac{16}{\pi}\sqrt{\frac{1+\sin (\alpha+y)}{1-\sin (\alpha+y)}}\left(|G_{1,0}|
+\frac{4\left|m_1B_{\alpha_1}^{-1}\right|}{\eta_{2}^{\alpha_{2}}(1-\sin (\alpha+y))^{\alpha_{2}}}|G_{2,0}|\right)e^{\eta_{2}t} \\
&\quad\times\left(\frac{L(\eta_{2}t\sin(\alpha))}{e^{2\pi(\pi/2-\alpha-\delta)/h_2}-1}+\frac{L(\eta_{2}t\sin(\alpha))}{e^{2\pi \alpha/h_2}-1}+\frac{1+L(\eta_{2}t\sin(\alpha))}{e^{\textcolor[rgb]{1.00,0.00,0.00}{\eta_{2}} t\sin{\alpha}\cosh(Nh_2)}}\right).
\end{aligned}\end{displaymath}

For $\left|G_2(t)-G^{\langle m\rangle}_{2,N}(t)\right|$, $m=1,2$, there are similar estimates hold.
\end{thm}

\begin{pf}
Let $v_{j}(t,\phi)$, $j=1,2$ be defined as in (\ref{exp:V1V2}) with $z(\phi)$ defined as (\ref{eq:paraboliccontour}) and (\ref{eq:hyperboliccontour}), respectively.

To proof the first part of the theorem, we need to perform error estimates about $DE^{\langle 1\rangle}_{+}$, $DE^{\langle 1\rangle}_{-}$, and $TE^{\langle 1\rangle}$.

Estimate of $DE^{\langle 1\rangle}_{+}$:
By (\ref{est:bound1}) in Lemma \ref{lemm:lemma3.1}, for $t>0$, $0<a<1$, there is
\begin{displaymath}\begin{aligned}
\int_{-\infty}^{\infty}|v_{1}(t,x+ia)|dx
&\leq \frac{\eta_{1}e^{\eta_{1}t}}{\pi}\left(\frac{16}{\eta_{1}(1-a)}|G_{1,0}|
  +\frac{64\left|m_1B_{\alpha_1}^{-1}\right|}{\eta_{1}^{1+\alpha_{2}}(1-a)^{1+2\alpha_{2}}}|G_{2,0}|\right)\int_{-\infty}^{\infty}e^{-\eta_{1}tx^{2}}dx\\
&\leq B\sqrt{\frac{\eta_{1}}{\pi}}t^{-1/2}e^{\eta_{1}t}\left(\eta_{1}^{-1}|G_{1,0}|
  +\eta_{1}^{-(1+\alpha_{2})}|G_{2,0}|\right)<\infty,
\end{aligned}\end{displaymath}
where $B$ is defined in (\ref{con:B}).
Then, according to Lemma \ref{lemm:lemma3.1}, we have
\begin{equation}\label{DE11}
DE^{\langle 1\rangle}_{+}\leq
B\sqrt{\frac{\eta_{1}}{\pi}t^{-1/2}e^{\eta_{1}t}}\left(\eta_{1}^{-1}|G_{1,0}|.
  +\eta_{1}^{-(1+\alpha_{2})}|G_{2,0}|\right)\frac{1}{e^{2\pi a/h_{1}}-1}.
\end{equation}
Thus,
\begin{equation}
DE^{\langle 1\rangle}_{+}=\mathscr{O}\left(e^{\eta_{1}t-2\pi a/h_{1}}\right),\quad h_{1}\rightarrow 0.
\end{equation}

Estimate of $DE^{\langle 1\rangle}_{-}$:
Similarly, for $t>0$, there holds
\begin{displaymath}\begin{aligned}
\int_{-\infty}^{\infty}|v_{1}(t,x-id)|dx
&\leq \frac{\eta_1e^{\eta_{1}(1+d)^{2}t}}{\pi}\left(16\eta_{1}^{-1}|G_{1,0}|
     +64\left|m_1B_{\alpha_1}^{-1}\right|\eta_{1}^{-(1+\alpha_{2})}|G_{2,0}|\right)\int_{-\infty}^{\infty}e^{-\eta_{1}tx^{2}}dx\\
&\leq B\sqrt{\frac{\eta_{1}}{\pi}}t^{-1/2}e^{\eta_{1}(1+d)^{2}t}\left(\eta_{1}^{-1}|G_{1,0}|+\eta_{1}^{-(1+\alpha_{2})}|G_{2,0}|\right)<\infty.
\end{aligned}\end{displaymath}
By Lemma \ref{lemm:lemma3.1}, we have
\begin{equation}\label{DE12}
DE^{\langle 1\rangle}_{-}\leq B\sqrt{\frac{\eta_{1}}{\pi}}t^{-1/2}e^{\eta_{1}(1+d)^{2}t}\left(\eta_{1}^{-1}|G_{1,0}|+\eta_{1}^{-(1+\alpha_{2})}|G_{2,0}|\right) \frac{1}{e^{2\pi d/h_{1}}-1},
\end{equation}
and
\begin{equation}\label{DE13}
DE^{\langle 1\rangle}_{-}=\mathscr{O}\left(e^{\eta_{1}(1+d)^{2}t-2\pi d/h_{1}}\right),\quad h_{1}\rightarrow 0.
\end{equation}

Denote $\omega(d)=\eta_{1}(1+d)^{2}t-2\pi d/h_{1}$. Then the `best' choice of $d$ is obtained by setting $\omega'(d)=0$ (see \cite{Weidemantrefe07}), which yields
\begin{equation}\label{eq:d}
d=\frac{\pi}{\eta_{1}th_{1}}-1.
\end{equation}
Now, we have
\begin{displaymath}
DE^{\langle 1\rangle}_{-}\leq B\sqrt{\frac{\eta_{1}}{\pi}}t^{-1/2}\left(\eta_{1}^{-1}|G_{1,0}|+\eta_{1}^{-(1+\alpha_{2})}|G_{2,0}|\right) \frac{e^{\pi^{2}/(\eta_{1}th_{1}^{2})}}{e^{2\pi^{2}/(\eta_{1}th_{1}^{2})-2\pi/h_1}-1}
\end{displaymath}
and
\begin{displaymath}
DE^{\langle 1\rangle}_{-}=\mathscr{O}\left(e^{-\pi^{2}/\left(\eta_{1}th_{1}^{2}\right)+2\pi/h_{1}}\right),\quad h_{1}\rightarrow 0.
\end{displaymath}

Estimate of $TE^{\langle 1\rangle}$: By the definition of $I_{1,h_{1}}$ and $I_{1,h_{1},N}$, we deduce
\begin{displaymath}
\left|I_{1,h_{1}}-I_{1,h_{1},N}\right|\leq h_{1}\sum_{k=N}^{\infty}(|v_{1}(t,kh_{1})|+|v_{1}(t,-kh_{1})|)
   \leq 2h_{1}\sum_{k=N}^{\infty}|v_{1}(t,kh_{1})|.
\end{displaymath}
According to (\ref{est:bound1}) in Lemma \ref{lemm:lemma3.1}, we have
\begin{align}
h_1\sum_{k=N}^{\infty}|v_{1}(t,kh_{1})|
&\leq \frac{\eta_{1}B e^{\eta_{1}t}}{\pi}\left(\eta_{1}^{-1}|G_{1,0}|+\eta_{1}^{-(1+\alpha_{2})}|G_{2,0}|\right)\int_{Nh_{1}}^{\infty}e^{-x^{2}\eta_{1}t}dx  \\ \nonumber
&\leq B\sqrt{\frac{\eta_{1}}{\pi}}t^{-1/2}\left(\eta_{1}^{-1}|G_{1,0}|+\eta_{1}^{-(1+\alpha_{2})}|G_{2,0}|\right)e^{\eta_{1}t\left(1-(Nh_{1})^{2}\right)}.
\end{align}
Thus,
\begin{equation}\label{TE11}
TE^{\langle 1\rangle}\leq B\sqrt{\frac{\eta_{1}}{\pi}}t^{-1/2}\left(\eta_{1}^{-1}|G_{1,0}|+\eta_{1}^{-(1+\alpha_{2})}|G_{2,0}|\right)e^{\eta_{1}t\left(1-(Nh_{1})^{2}\right)}
\end{equation}
and
\begin{equation}\label{TEE1}
TE^{\langle 1\rangle}=\mathscr{O}\left(e^{\eta_{1}t\left(1-(Nh_{1})^{2}\right)}\right),\quad N\rightarrow +\infty.
\end{equation}
Combining (\ref{DE11}), (\ref{DE12}), and (\ref{TE11}) results in the first part of the theorem.

As for the hyperbolic integral contour with the strip $S^{\langle2\rangle}$, similar to the previous analyses, we will directly give the corresponding results in the sequence.

Estimate of $DE^{\langle 2\rangle}_{+}$:
By (\ref{est:bound2}) in Lemma \ref{lemm:lemma3.1} and Lemma \ref{lem:lemma2}, for $t>0$ and $0<y<\pi/2-\alpha-\delta$, there holds
\begin{equation}\label{DE21}
DE^{\langle 2\rangle}_{+}
\leq  \frac{16}{\pi}\sqrt{\frac{1+\sin (\alpha+y)}{1-\sin (\alpha+y)}}\left(|G_{1,0}|
     +\frac{4\left|m_1B_{\alpha_1}^{-1}\right|}{\eta_{2}^{\alpha_{2}}(1-\sin (\alpha+y))^{\alpha_{2}}}|G_{2,0}|\right)\frac{e^{\eta_{2}t}L(\eta_{2}t\sin(\alpha))}{e^{2\pi(\pi/2-\alpha-\delta)/h_{2}}-1},
\end{equation}
and
\begin{equation}
DE^{\langle 2\rangle}_{+}=\mathscr{O}\left(e^{\eta_{2}t-2\pi(\pi/2-\alpha-\delta)/h_{2}}\right),\quad h_{2}\rightarrow0.
\end{equation}

Estimate of $DE^{\langle 2\rangle}_{-}$:
Similarly, there are
\begin{equation}\label{DE22}
DE^{\langle 2\rangle}_{-}
\leq  \frac{16}{\pi}\sqrt{\frac{1+\sin (\alpha+y)}{1-\sin (\alpha+y)}}\left(|G_{1,0}|
    +\frac{4\left|m_1B_{\alpha_1}^{-1}\right|}{\eta_{2}^{\alpha_{2}}(1-\sin (\alpha+y))^{\alpha_{2}}}|G_{2,0}|\right)\frac{e^{\eta_{2}t}L(\eta_{2}t\sin(\alpha))}{e^{2\pi\alpha/h_{2}}-1},
\end{equation}
and
\begin{equation}
DE^{\langle 2\rangle}_{-}=\mathscr{O}\left(e^{\eta_{2}t-2\pi\alpha/h_{2}}\right),\quad h_{2}\rightarrow0.
\end{equation}

Estimate of $TE^{\langle 2\rangle}$:
From (\ref{est:bound2}) and Lemma \ref{et:lemma}, there holds
\begin{displaymath}\begin{aligned}
 h_{2}\sum_{N}^{\infty}|v_{1}(t,kh_{2})|
& \leq \frac{8}{\pi}\sqrt{\frac{1+\sin (\alpha+y)}{1-\sin (\alpha+y)}}\left(|G_{1,0}|
     +\frac{4\left|m_1B_{\alpha_1}^{-1}\right|}{\eta_{2}^{\alpha_{2}}(1-\sin (\alpha+y))^{\alpha_{2}}}|G_{2,0}|\right)\int_{Nh_{2}}^{+\infty}e^{\eta_{2}t-\eta_{2}t\sin\alpha\cosh x}dx\\
& \leq \frac{8}{\pi}(1+L(\eta_{2}t\sin\alpha))\sqrt{\frac{1+\sin (\alpha+y)}{1-\sin (\alpha+y)}}\left(|G_{1,0}|
     +\frac{4\left|m_1B_{\alpha_1}^{-1}\right|}{\eta_{2}^{\alpha_{2}}(1-\sin (\alpha+y))^{\alpha_{2}}}|G_{2,0}|\right)e^{\eta_{2}t-\eta_{2}t\sin(\alpha)\cosh(Nh_{2})}.
\end{aligned}\end{displaymath}
Thus
\begin{equation}\label{TEE2}
TE^{\langle 2\rangle} \leq \frac{16}{\pi}(1+L(\eta_{2}t\sin\alpha))\sqrt{\frac{1+\sin (\alpha+y)}{1-\sin (\alpha+y)}}\left(|G_{1,0}|
+\frac{4\left|m_1B_{\alpha_1}^{-1}\right|}{\eta_{2}^{\alpha_{2}}(1-\sin (\alpha+y))^{\alpha_{2}}}|G_{2,0}|\right)e^{\eta_{2}t-\eta_{2}t\sin(\alpha)\cosh(Nh_{2})},
\end{equation}
and
\begin{equation}
TE^{\langle 2\rangle}=\mathscr{O}\left(e^{\eta_{2}t(1-\sin{\alpha}\cosh(h_{2}N))}\right),\quad N\rightarrow +\infty.
\end{equation}
Together with (\ref{DE21}), (\ref{DE22}) and (\ref{TEE2}), we finish the second part of the theorem.

Similar estimates on $\left|G_2(t)-G^{\langle m\rangle}_{2,N}(t)\right|$ can be obtained.
\end{pf}

Next, we determine the optimal parameters in the integral contours $\Gamma_1$ and $\Gamma_2$. Reference \cite{Weidemantrefe07} has provided a technical method to determine these parameters under ideal conditions, i.e., asymptotically balancing $DE^{\langle m\rangle}_{+}$, $DE^{\langle m\rangle}_{-}$, and $TE^{\langle m\rangle}$, $m=1, 2$. Based on this idea, we will optimize these parameters in our situations.

For $\Gamma_1$, by asymptotically balancing $DE^{\langle 1\rangle}_{+}$, $DE^{\langle 1\rangle}_{-}$, and $TE^{\langle 1\rangle}$, it needs
\begin{equation}\label{tpeq:err}
\eta_{1}t-\frac{2\pi a}{h_{1}}=\frac{-\pi^{2}}{\eta_{1}th_{1}^{2}}+\frac{2\pi}{h_{1}}=\eta_{1}t\left(1-(Nh_{1})^{2}\right).
\end{equation}
Then we obtain the optimal parameters used in the contour $\Gamma_1$, i.e.,
\begin{equation}\label{optm:param1}
\eta_{1}^{\ast}= \frac{\pi \sqrt{2aq^{3}}}{2a}\frac{N}{t},~~{\rm and}~ h_{1}^{\ast}=\frac{\sqrt{2aq}}{q}\frac{1}{N},
\end{equation}
where $q:=1+a-\sqrt{a^{2}+2a}$ with $1/4<a<1$.
With these optimal parameters, the corresponding convergence order of the CIM with the parabolic contour $\Gamma_{1}$ are
\begin{displaymath}
E^{\langle 1\rangle}_{j,N}=\mathscr{O}\left(e^{-\left(\pi\sqrt{2aq}-\pi\sqrt{2aq^{3}}/(2a)\right)N}\right),~ j=1,2.
\end{displaymath}

Furthermore, as mentioned in \cite{Weidemantrefe07}, the parameter $\eta_{1}^{\ast}$ in (\ref{optm:param1}) depends on time $t$, which means that the integral contour changes with time $t$. The ideal situation is that we find a fixed integral contour that satisfies the condition (\ref{condition02}) which  does not change over time $t$. Reviewing the error estimates $DE^{\langle 1\rangle}_{+}$, $DE^{\langle 1\rangle}_{-}$, and $TE^{\langle 1\rangle}$, it can be found that $DE^{\langle 1\rangle}_{+}$ and $DE^{\langle 1\rangle}_{-}$ increase with $t$, and $TE$ decreases with $t$. If we want a small absolute error on the interval $t_{0}\leq t\leq t_{1}=T$, $t_{0}>0$, we can modify (\ref{tpeq:err}) as
\begin{equation}\label{tpeq1:err}
\eta_{1}t_{1}-\frac{2\pi a}{h_{1}}=\frac{-\pi^{2}}{\eta_{1}t_{1}h_{1}^{2}}+\frac{2\pi}{h_{1}}=\eta_{1}t_{0}\left(1-(Nh_{1})^{2}\right).
\end{equation}
Denote $\Lambda=t_{1}/t_{0}$, which increases from $1$. After solving (\ref{tpeq1:err}), we have
\begin{equation}\label{eq:eta_1}
\eta_{1}=\frac{\pi q\sqrt{q^{2}(1-\Lambda)+2a\Lambda q}}{q(1-\Lambda)+2a\Lambda }\frac{N}{t_{1}},\quad h_{1}=\frac{\sqrt{q^{2}(1-\Lambda)+2a\Lambda q}}{q}\frac{1}{N},
\end{equation}
and the corresponding convergence order of the CIMs are
\begin{equation}
E^{\langle 1\rangle}_{j,N}=\mathscr{O}\left(e^{-P(\Lambda)N}\right),~ ~j=1,2,~ N\rightarrow +\infty,
\end{equation}
where $P(\Lambda)=\frac{\pi(q-2a)\sqrt{q^{2}(1-\Lambda)+2a\Lambda q}}{q(\Lambda-1)-2a\Lambda}$.

For $\Gamma_2$, similar to the previous analyses of $\Gamma_1$, when $t\in[t_{0}, t_{0}\Lambda]$, $t_{0}>0, T=t_{1}=t_{0}\Lambda$, the discretization error $DE^{\langle 2\rangle}_{-}$ increases with $t$ and the truncation error $TE^{\langle 2\rangle}$ decreases with $t$. Thus, $DE^{\langle 2\rangle}_{-}$ and $TE^{\langle 2\rangle}$ can be modified as
\begin{displaymath}
DE_{-}^{\langle 2\rangle ~*}=\mathscr{O}\left(e^{\eta_{2}t_{1}-2\pi\alpha/h_{2}}\right),
\quad TE^{\langle 2\rangle ~*}=\mathscr{O}\left(e^{\eta_{2}t_{0}(1-\sin{\alpha}\cosh(h_{2}N))}\right).
\end{displaymath}
By asymptotically balancing $DE^{\langle 2\rangle}_{+}$, $DE_{-}^{\langle 2\rangle ~*}$ and $TE^{\langle 2\rangle ~*}$, there holds
\begin{displaymath}
\frac{-2\pi(\pi/2-\alpha-\delta)}{h_{2}}= \eta_{2}t_{1}-\frac{2\pi\alpha}{h_{2}}=\eta_{2}t_{0}(1-\sin{\alpha}\cosh(h_{2}N)).
\end{displaymath}
Solving it results in
\begin{equation}\label{eq:eta_2}
h_2=\frac{A(\alpha)}{N},~~\eta_{2}=\frac{4\pi\alpha-\pi^{2}+2\pi\delta}{A(\alpha)}\frac{N}{t_{1}}, ~ ~ {\rm and}~ A(\alpha)=\cosh^{-1}\left(\frac{(\pi-2\alpha-2\delta)\Lambda+(4\alpha-\pi+2\delta)}{(4\alpha-\pi+2\delta)\sin(\alpha)}\right),
\end{equation}
with  and
\begin{displaymath}
E_{j,N}^{\langle 2\rangle}=\mathscr{O}\left(e^{-Q(\alpha)N}\right),~ j=1,2,~ N\rightarrow +\infty,
\end{displaymath}
where $Q(\alpha)=\frac{\pi^{2}-2\pi\alpha-2\pi\delta}{A(\alpha)}$. For fixed $\Lambda$ and $\delta$, the optimal parameter $\alpha$ can be obtained by maximizing $Q(\alpha)$, which is similar to the results of \cite{Weidemantrefe07}.

To sum up: The specific strategy we determine the integral contours $\Gamma_1$ is as follows
\begin{enumerate}[label=\alph*)]
    \item given $a$ with $1/4<a<1$; compute $q=1+a-\sqrt{a^{2}+2a}$;
    \item determine $\eta_1$ and $h_1$ by (\ref{eq:eta_1});
    \item let $\phi=x+iy$ with $y=0$ and $x\in \mathbb{R}$;
    \item determine $\Gamma_1$ as defined in (\ref{eq:paraboliccontour}).
\end{enumerate}
The specific strategy we determine the integral contours $\Gamma_2$ is as below
\begin{enumerate}[label=\alph*)]
    \item given an interval $[t_0,t_1]$; compute $\Lambda=t_1/t_0$;
    \item given $\delta$ with $0<\delta<\pi/2$; determine the optimal parameter $\alpha^*$ by maximizing $Q(\alpha)=\frac{\pi^{2}-2\pi\alpha-2\pi\delta}{A(\alpha)}$;
    \item compute $A(\alpha^*)$ and then $h_2$ and $\eta_2$ in (\ref{eq:eta_2});
    \item let $\phi=x+iy$ with $x\in \mathbb{R}$ and a given $y$, $0<y<\min\{\alpha, \pi/2-\alpha-\delta\}$;
    \item determine $\Gamma_2$ as defined in (\ref{eq:hyperboliccontour}).
\end{enumerate}

Through the above analyses, it can be found that the CIMs constructed in this paper with given optimal step-sizes and parameters have the convergence order of $\mathscr{O}(e^{-cN})$. That is, the CIMs in our paper have spectral accuracy.

\subsection{The time-matching schemes}
\label{subsec:PC}
In this subsection, in order to verify the high numerical performance of the CIMs (\ref{eq:computG1}) with the determined integral contours $\Gamma_1$ and $\Gamma_2$, we also show another numerical methods to solve (\ref{eq:DIVFeynman-kac}), i.e., the time-marching schemes (TMs) e.g. \cite{Wang95}, etc.

%

Let $h$ be the discrete stepsize and $t_{n}=nh$, $n=0,1,2,\cdot\cdot\cdot,M$, $M=\frac{T}{h}$. Integrating both sides of (\ref{bianxing}) from $0$ to $t$ and letting $t=t_{n+1}$, we get the integral form (\ref{sys:integral}) and (\ref{sys:integral2}).


The specific TMs  are shown in \ref{sec:Pre-Corre}, in which the following two formulas are used.
\begin{equation}\begin{aligned}\label{eq:integral}
  \int_{0}^{t}\mathfrak{D}_{s}^{1-\alpha_{j}}G_{j}(.,\rho )ds
= &\frac{1}{\Gamma(\alpha_{j})}\int_{0}^{t}\frac{\partial}{\partial s}\int_{0}^{s}\frac{e^{-(s-\tau)\rho U(j)}}{(s-\tau)^{1-\alpha_{j}}}G_{j}(\rho ,\tau)d\tau ds
  +\frac{\rho U(j)}{\Gamma(\alpha_{j})}\int_{0}^{t}\int_{0}^{s}\frac{e^{-(s-\tau)\rho U(j)}}{(s-\tau)^{1-\alpha_{j}}}G_{j}(\rho ,\tau)d\tau ds\\
= &\frac{1}{\Gamma(\alpha_{j})}\left(\int_{0}^{t}\frac{e^{-(t-\tau)\rho U(j)}}{(t-\tau)^{1-\alpha_{j}}}G_{j}(\rho ,\tau)d\tau
+(\rho U(j))^{1-\alpha_{j}}\gamma(\alpha_{j},1)\int_{0}^{t}G_{j}(\rho ,\tau)d\tau\right),~ j=1,2,
\end{aligned}\end{equation}
where $\gamma(\alpha_{j},1)$ are the incomplete gamma function with different parameters $\alpha_{j}, j=1,2$. For the  first term of the right hand of \ref{eq:integral}, we perform the linear interpolation on the integrand, that is
\begin{displaymath}\begin{aligned}
\int_{t_{n}}^{t_{n+1}}e^{-\rho U(j)(t_{n+1}-\tau)}(t_{n+1}-\tau)^{\alpha_{j}-1}G_{j}(\rho ,\tau)d\tau
&\approx\int_{t_{n}}^{t_{n+1}}(t_{n+1}-\tau)^{\alpha_{j}-1}\frac{G_{j}(.,t_{n+1})(\tau-t_{n})+e^{-\rho U(j)h}G_{j}(.,t_{n})(t_{n+1}-\tau)}{h}d\tau\\
&=\frac{h^{\alpha_{j}}}{\alpha_{j}(\alpha_{j}+1)}\left(\alpha_{j}e^{-\rho U(j)h}G_{j}(.,t_{n})+G_{j}(.,t_{n+1})\right),~~ j=1,2.
\end{aligned}\end{displaymath}

Based on these results, after doing some calculations, the time-marching schemes of (\ref{eq:DIVFeynman-kac}) are designed as in (\ref{for:n0}), (\ref{for:nn}).

\section{Numerical Results}
\label{sec:NumericalResults}
In this section, we use two examples to evaluate the effectiveness of the CIM.  We choose $p=2/3$, $b=3/4$, $B_{\alpha_1}^{-1}=B_{\alpha_2}^{-1}=1$, and the optimal parameters are given in (\ref{eq:eta_1}) and (\ref{eq:eta_2}).
Take $a=0.9875$ and $\delta=0.1123$ such that the contours $\Gamma_1$ and $\Gamma_2$ satisfy the condition (\ref{condition02}).
Here, we remark that the computing environment for all examples is \emph{Intel(R) Core(TM) i7-7700 CPU @3.60GHz, MATLAB R2018a}.

\subsection{Example 1}
\label{exam:01}
\begin{figure}[h]
\centering
\begin{minipage}[c]{0.32\textwidth}
 \centering
 \centerline{\includegraphics[height=4.2cm,width=5.7cm]{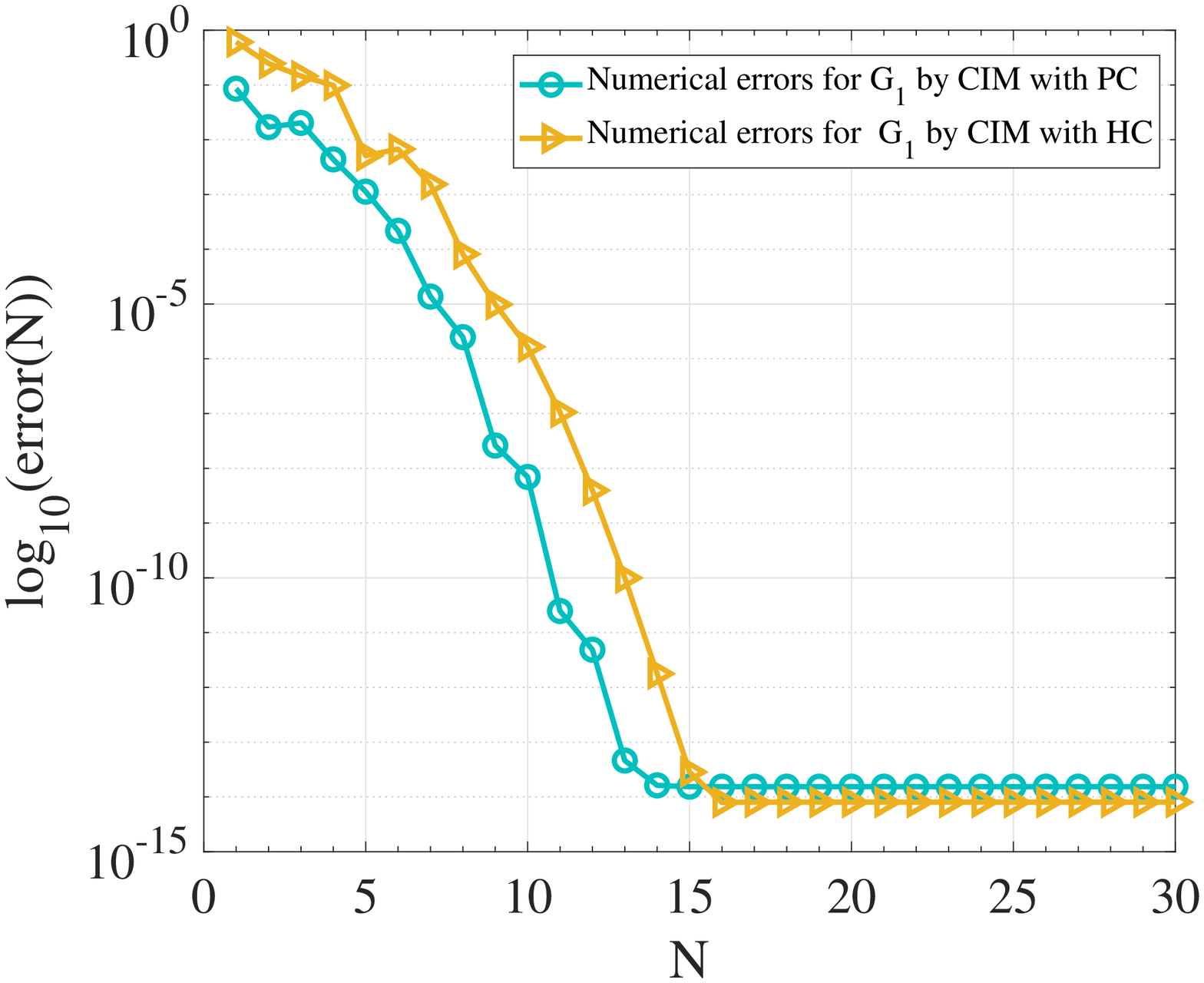}}
\end{minipage}
\begin{minipage}[c]{0.32\textwidth}
 \centering
 \centerline{\includegraphics[height=4.2cm,width=5.7cm]{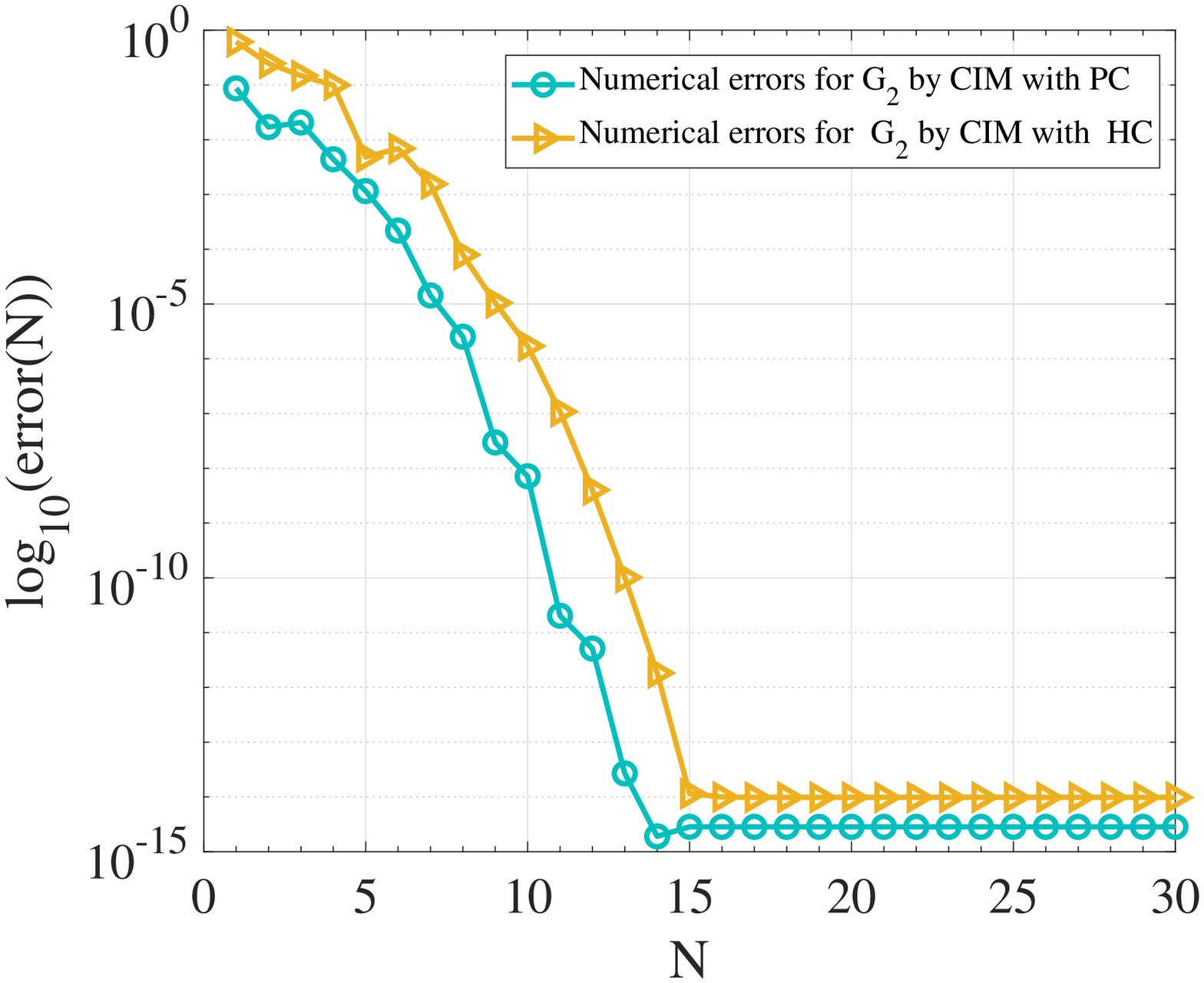}}
\end{minipage}
\begin{minipage}[c]{0.32\textwidth}
 \centering
 \centerline{\includegraphics[height=4.2cm,width=5.7cm]{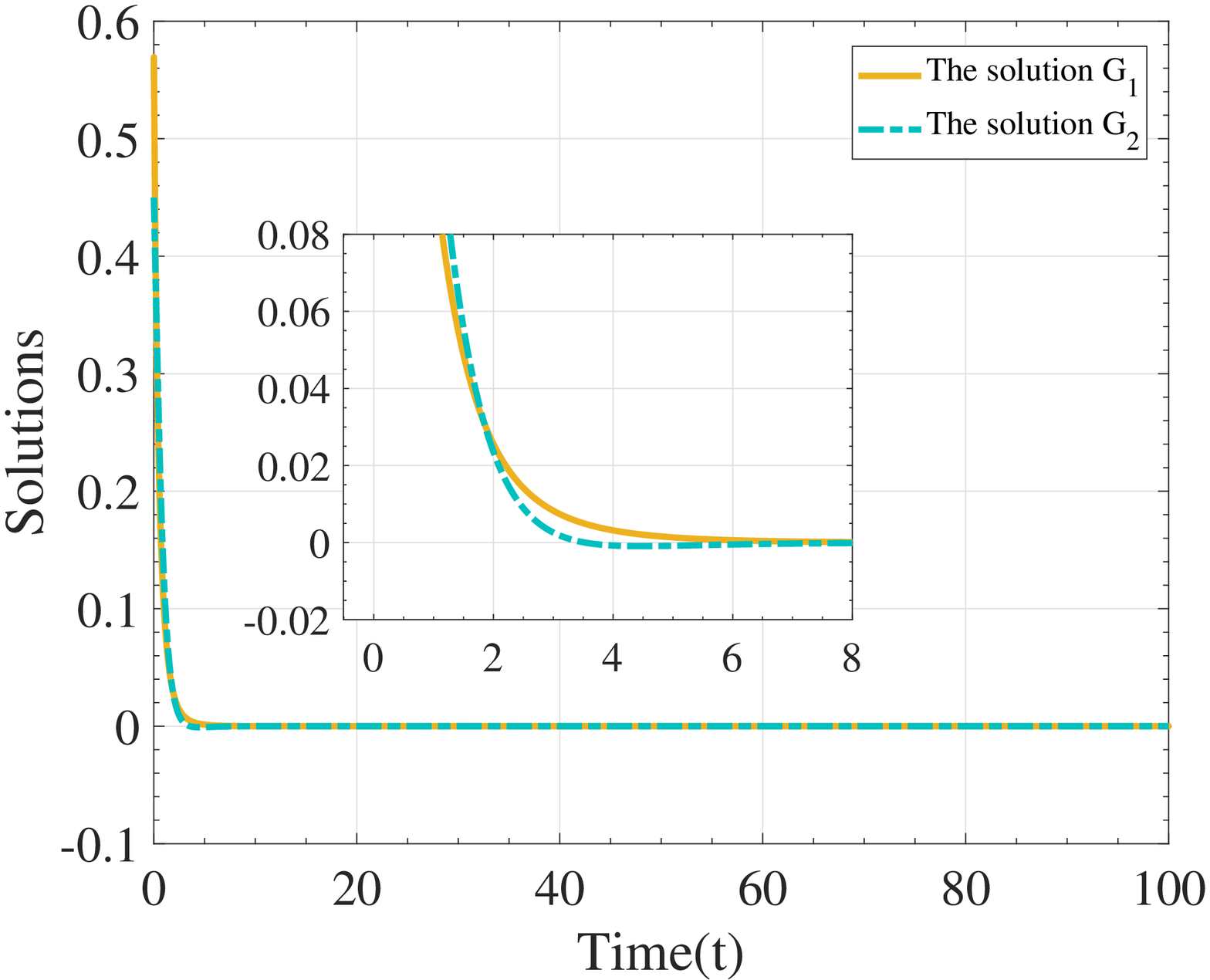}}
\end{minipage}
\caption{The absolute errors at time $t=100$ and the numerical solutions of the system (\ref{eq:DIVFeynman-kac}); (\emph{left}) and (\emph{center}) are the absolute errors for the CIM with $\Gamma_1$ and $\Gamma_2$, respectively; The parameters are set as $N=30$, $\Lambda$=5, $\alpha_{1}=0.6$, $\alpha_{2}=0.4$, $U(1)=1$, $U(2)=1$, $\rho=1.5$, and the initial values $G_{1,0}=0.55$, $G_{2,0}=0.45$; (\emph{right}) are the numerical solutions obtained by using TMs with the same parameters and the number of discrete points $M=2^{15}$ and $T=100$.}
\label{Fig:example11}
\end{figure}

\begin{figure}[h]
\centering
\begin{minipage}[c]{0.32\textwidth}
 \centering
 \centerline{\includegraphics[height=4.2cm,width=5.7cm]{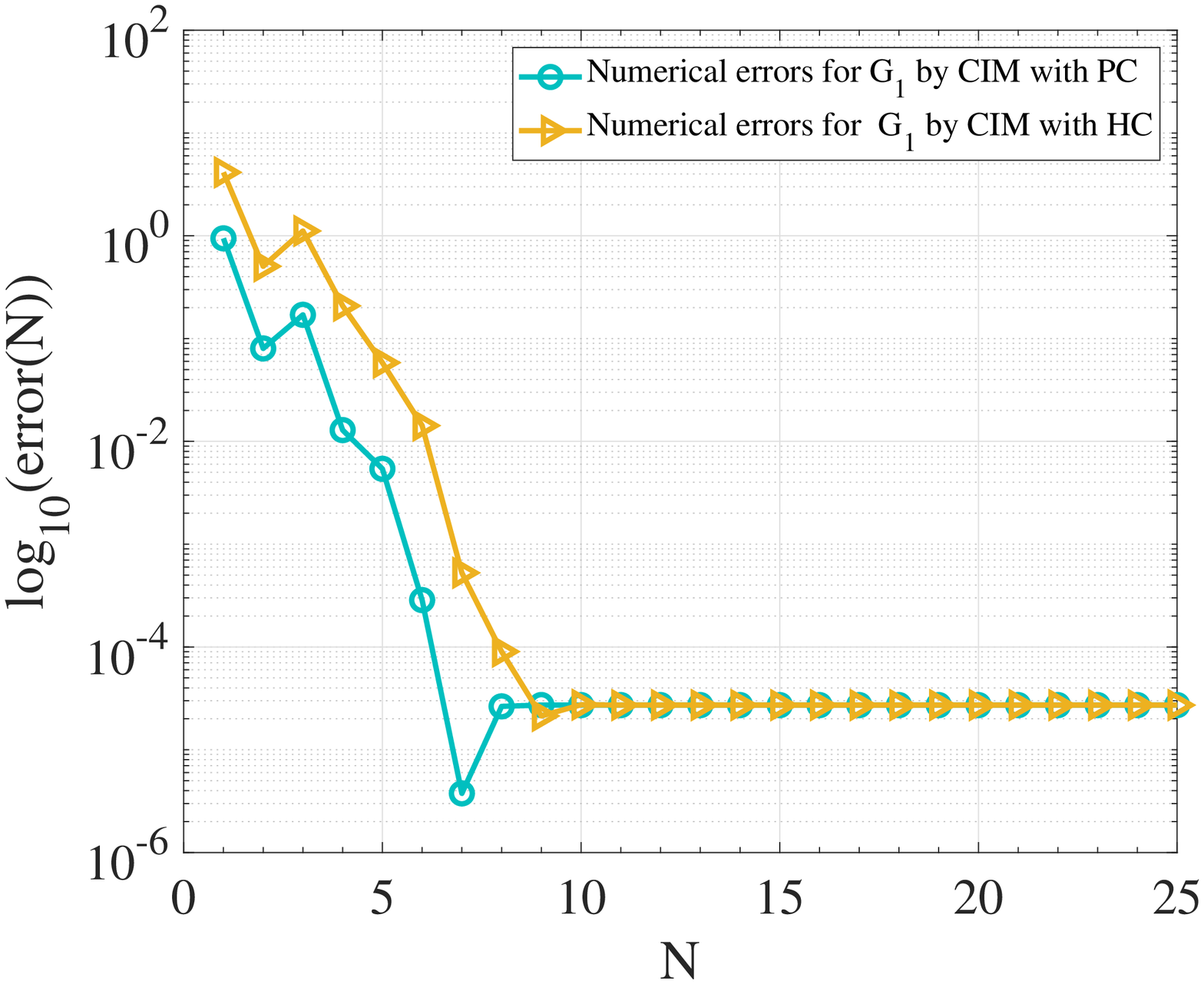}}
\end{minipage}
\begin{minipage}[c]{0.32\textwidth}
 \centering
 \centerline{\includegraphics[height=4.2cm,width=5.7cm]{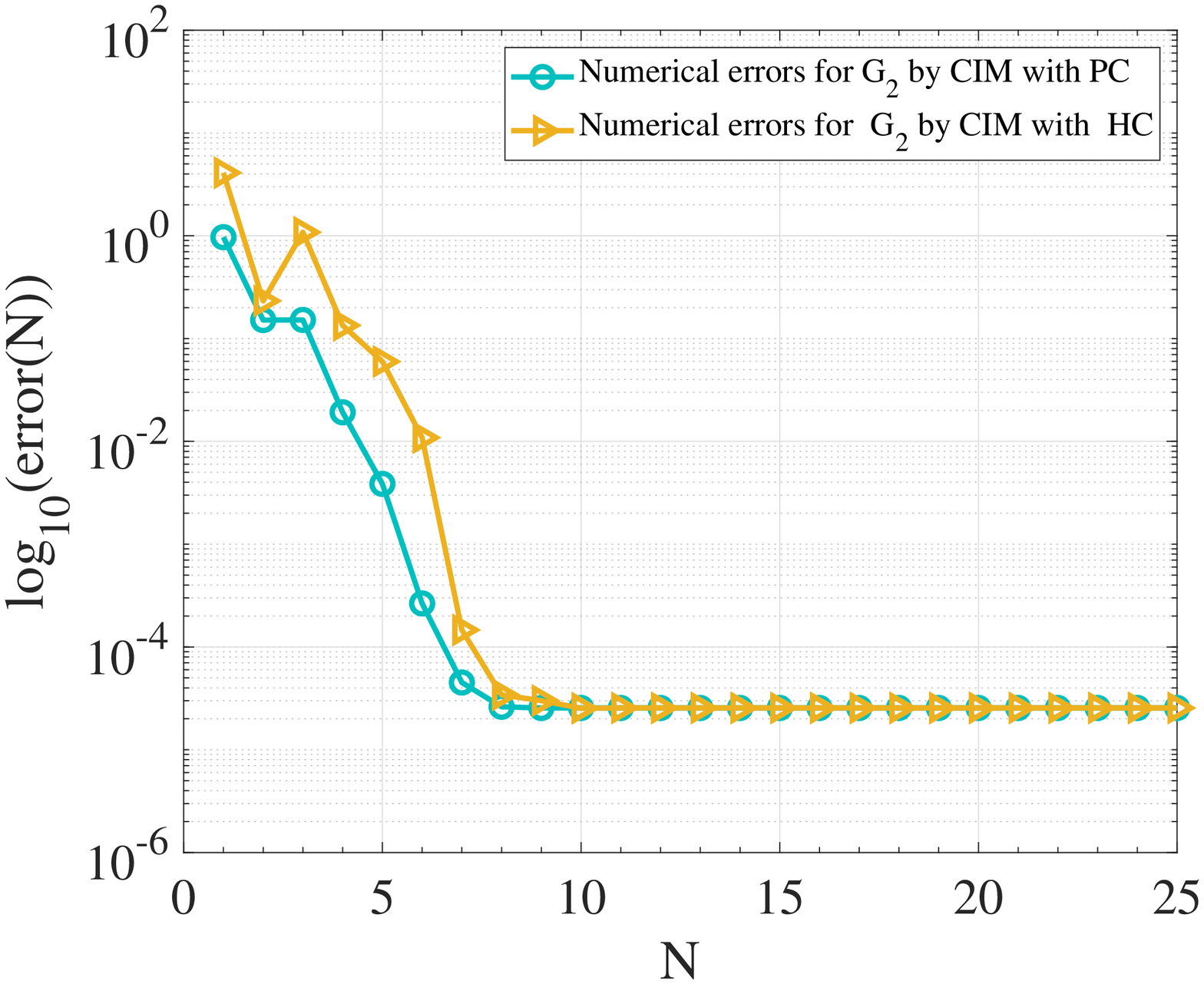}}
\end{minipage}
\begin{minipage}[c]{0.32\textwidth}
 \centering
 \centerline{\includegraphics[height=4.2cm,width=5.7cm]{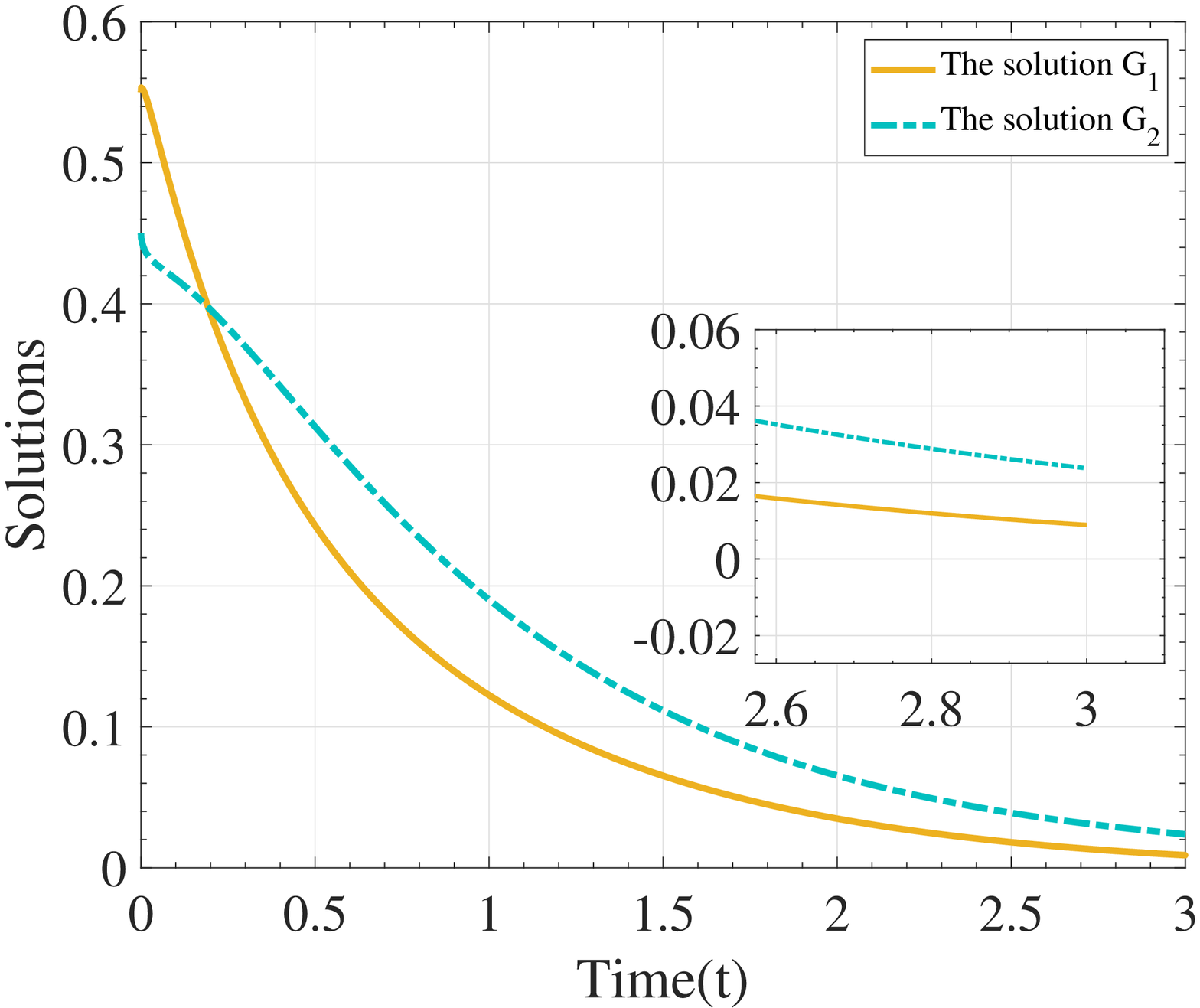}}
\end{minipage}
\caption{The absolute errors at time $t=3$ and the numerical solutions of the system (\ref{eq:DIVFeynman-kac}); (\emph{left}) and (\emph{center}) are the absolute errors of the CIM with $\Gamma_1$ and $\Gamma_2$, respectively; The parameters are set as $N=25$, $\Lambda$=5, $\alpha_{1}=0.82$, $\alpha_{2}=0.59$, $U(1)=0.89$, $U(2)=0.68$, $\rho=1.5$, and the initial values $G_{1,0}=0.55$ and $G_{2,0}=0.45$; (\emph{right}) are the numerical solutions obtained by using TMs with the same parameters and the number of discrete points $M=2^{12}$ and $T=3$.}
\label{Fig:example}
\end{figure}

This example is used to verify that the CIMs (\ref{eq:computG1}) have spectral accuracy. We choose the absolute error as a function of $N$, i.e.,
\begin{equation}
  error(N)=\max\limits_{t_0\leq t\leq \Lambda t_0}\left|G_{j}(.,t)-G^{\langle m \rangle}_{j,N}(.,t)\right|,~ ~ m, j=1,2,
\end{equation}
where $G_{j}(.,t)$ is the reference solution computed by the time-marching schemes with much small stepsize, and $G^{\langle m \rangle}_{j,N}(.,t)$ is the numerical solution obtained by the CIM with parabolic contour $\Gamma_1$ and hyperbolic contour $\Gamma_2$, respectively. Then, the absolute errors of the CIMs at different given times and the numerical solutions of the system are shown in Figure \ref{Fig:example11} and Figure \ref{Fig:example}.

Given a specific accuracy to be reached, the number of discrete points and the time cost of the CIMs and TMs are shown in Table \ref{Tab:01}, in which the parameters used are consistent with those in Figure 3, and we treat the numerical solutions obtained by TMs with $M=2^{12}$.


\begin{table*}[t!]
\scriptsize
\caption{The CPU time cost of CIM and TMS for the system (\ref{eq:DIVFeynman-kac}) at time $t=3$.}
\begin{center}
\begin{tabular}{cccccccc}
\toprule
\multirow{2}{*}{Accuracy}&\multirow{2}{*}{PDF}&\multicolumn{2}{c}{CIM-PC}&\multicolumn{2}{c}{CIM-HC}&\multicolumn{2}{c}{TMs}\\
                                               \cmidrule(r){3-4}            \cmidrule(r){5-6}            \cmidrule(r){7-8}
                                &                &N      &time(s)        &N      &time(s)        &M        &time (s)\\
\midrule
\multirow{2}{*}{$10^{-2}$}      &$G_{1}(.,t)$    &5      &4.1903e-03     &3      &1.0031e-02     &10       &2.0175e-03\\
                                &$G_{2}(.,t)$    &5      &3.4310e-03     &3      &1.2785e-02     &24       &5.2115e-03\\
\midrule
\multirow{2}{*}{$10^{-3}$}      &$G_{1}(.,t)$    &9      &1.2757e-03     &7      &1.6547e-03     &10       &2.2867e-03\\
                                &$G_{2}(.,t)$    &9      &1.1832e-03     &7      &2.0558e-03     &283      &1.6356e-01\\
\midrule
\multirow{2}{*}{$10^{-4}$}      &$G_{1}(.,t)$    &13     &3.0482e-03     &8      &1.2326e-03     &630      &6.7298e-01\\
                                &$G_{2}(.,t)$    &14     &2.3619e-03     &8      &1.4250e-03     &1953     &5.9451e+00\\
\bottomrule
\end{tabular}
\end{center}
\label{Tab:01}
\end{table*}

One can see from Figures \ref{Fig:example11}-\ref{Fig:example3} as well as Table \ref{Tab:01} that the CIMs with parabolic contour and hyperbolic contour are much effective and time-saving in solving the Feynman-Kac equation with two internal states. Also they have spectral accuracy.

\subsection{Example 2}
\begin{figure}[h]
\centering
\begin{minipage}[c]{0.40\textwidth}
 \centering
 \centerline{\includegraphics[height=5.0cm,width=6.5cm]{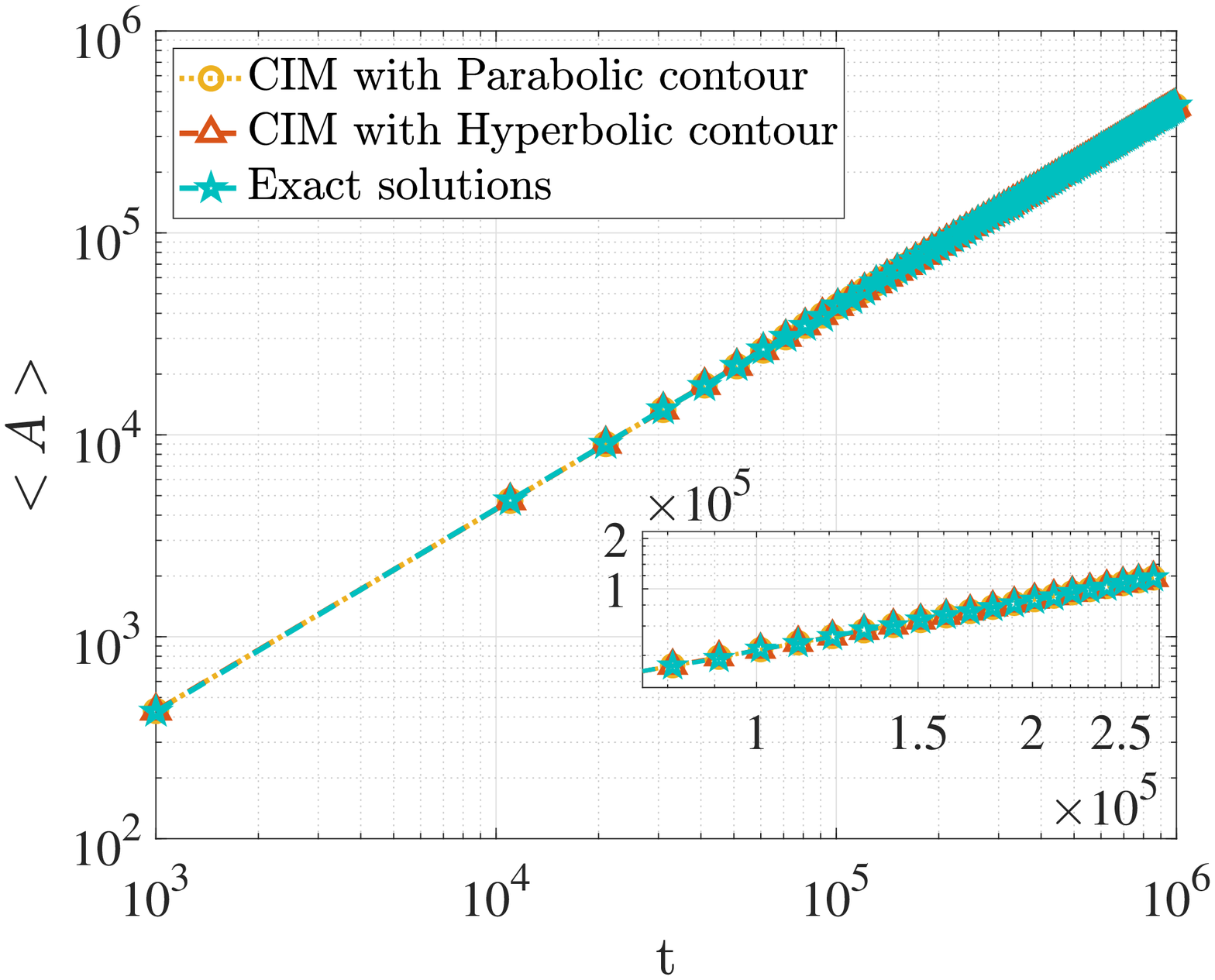}}
\end{minipage}
\begin{minipage}[c]{0.40\textwidth}
 \centering
 \centerline{\includegraphics[height=5.0cm,width=6.5cm]{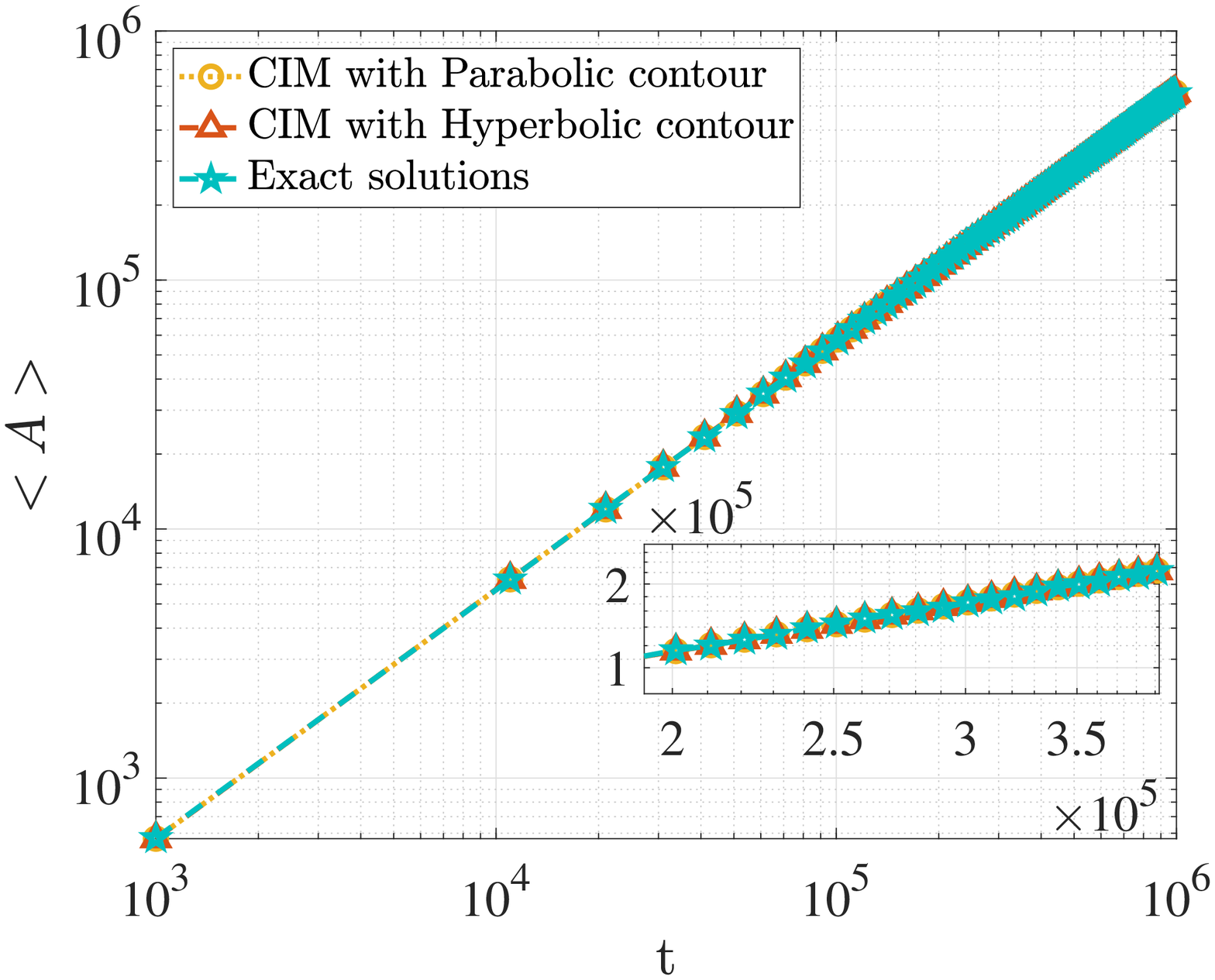}}
\end{minipage}
\caption{The simulations of $\langle A\rangle$ (in log-log scale); (\emph{left})  $\langle A\rangle$ for the first internal state, i.e., $U(1)=1$, $U(2)=0$; (\emph{right})  $\langle A\rangle$ for the second internal state, i.e., $U(1)=0$, $U(2)=1$. For both (\emph{left}) and (\emph{right}), we take $N=100$, $\Lambda=50$, $\alpha_{1}=\alpha_{2}=0.5$, $G_{1,0}=0.45$, $G_{2,0}=0.55$, $\varepsilon_{1}=1/2$, $\varepsilon_{2}=1/2$, $t\in[1E+3,1E+6]$.}
\label{Fig:example3}
\end{figure}

As a physical application, we calculate the average occupation time of each internal state by using the CIMs. As the matter of fact, the average occupation time of the first internal state can be calculated as the solution of the system (\ref{eq:Feynmankac}) by taking
\begin{equation}\label{bigu}
U[j(\tau)]=
\left\{
\begin{aligned}
&1,~ j(\tau)=1,\\
&0,~ {\rm else},
\end{aligned}\right.
\end{equation}
in the functional $A=\int_{0}^{\tau}[U(j(\tau))]d\tau$. While taking $U(1)=0$ and $U(2)=1$, it works for the second internal state.
Take $\alpha_{1}=\alpha_{2}$. According to the theoretical results presented in \cite{xu18}, the average occupation time of the first state is
\begin{equation}\label{eq:As}
\langle A\rangle\sim\frac{\varepsilon_{1}}{\varepsilon_{1}+\varepsilon_{2}}t,\quad~ ~ {\rm for}~ ~ {\rm large}~ ~ t,
\end{equation}
where $\varepsilon_{1}$ and $\varepsilon_{2}$ are the initial distributions (replacing the coefficient by $\varepsilon_{2}/(\varepsilon_{1}+\varepsilon_{2})$ for the second internal state). In this paper, $\langle A\rangle$ calculated by using the fact
\begin{equation}\label{eq:NAs2}
\langle A\rangle=-\frac{\partial}{\partial\rho} g(\rho,t)\Big|_{\rho=0}
\end{equation}
with $g(\rho,t)=G_{1}(\rho,t)+G_{2}(\rho,t)$. More specifically, after differentiating (\ref{eq:solutionsinlaplace}) w.r.t. $\rho$ and setting $\rho=0$, we then solve the system by the CIMs. At this time, one can check that the analytical domain of the new system is no smaller than the original one, therefore the discussions in above sections still work. See the simulation results in Figure \ref{Fig:example3}, which further verifies the theoretical predictions.

\section{Conclusions}
\label{sec:conclusions}
The Feynman-Kac equation with two internal states describes the functional distribution of the particle's internal states.  This paper presents the regularity analyses for the system, and built the CIMs with numerical stability analysis and error estimates. With the reference solutions provided by the time-marching schemes, the performances (of the CIMs) on spectral accuracy, low computational complexity, and small memory requirement, etc, are obtained. As one of the physical applications, by using the CIMs, we calculate the average occupation time of the first and second state of the stochastic process with two internal states.

\appendix
\section{The proofs of the results presented in Section \ref{subsec:Regularity}}
\label{subsec:RegularityProof}
The techniques of the proofs 
are inspired by \cite{Dengli18}.
\label{sec:Appendix}

\subsection{The proof of Lemma \ref{lemm:lemma2.4}}

\begin{pf}
For $z\in\Sigma_{\theta,\delta}$, there is
\begin{displaymath}
\frac{\left|(z+\rho U(j))^{-\alpha_{j}}\right|}{|z|^{-\alpha_{j}}}=\left(\frac{|z|}{|z+\rho U(j)|}\right)^{\alpha_{j}}\leq
      \left(\frac{|z|}{|z|-|\rho U(j)|}\right)^{\alpha_{j}}.
\end{displaymath}
Since $\delta>2|\rho_{s}U(j)|$, it has
\begin{displaymath}
|z|-|\rho U(j)|\geq|z|-\frac{1}{2}|z|=\frac{1}{2}|z|.
\end{displaymath}
Then
\begin{displaymath}
\frac{|(z+\rho U(j))^{-\alpha_{j}}|}{|z|^{-\alpha_{j}}}\leq(2)^{\alpha_{j}}\leq2.
\end{displaymath}
\end{pf}
\subsection{The proof of Lemma \ref{lemm:lemma2.5}}

\begin{pf}
For $z\in\Sigma_{\theta,\delta}$, let
\begin{displaymath}
\left((z+\rho U(j))^{\alpha_{j}}-m_{j}B_{\alpha_{j}}^{-1}\right)u_j=v_j,~  ~j=1,2.
\end{displaymath}
Then
\begin{displaymath}
u_j=\left((z+\rho U(j))^{\alpha_{j}}\right)^{-1}v_j+m_{j}B_{\alpha_{j}}^{-1}\left((z+\rho U(j))^{\alpha_{j}}\right)^{-1}u_j, ~  ~ j=1,2.
\end{displaymath}
Taking modulus on both sides of the above equality and by Lemma \ref{lemm:lemma2.4}, we have
\begin{displaymath}\begin{aligned}
\left|u_j\right|
 \leq \left|\left((z+\rho U(j))^{\alpha_{j}}\right)^{-1}\right|\left|v_j\right|
+\left|m_{j}B_{\alpha_{j}}^{-1}\left((z+\rho U(j))^{\alpha_{j}}\right)^{-1}\right| \left|u_j\right|
 \leq 2 |z|^{-\alpha_{j}}\left|v_j\right|+ 2 \left|m_{j}B_{\alpha_{j}}^{-1}\right||z|^{-\alpha_{j}}\left|u_j\right|, j=1,2.
\end{aligned}\end{displaymath}
According to the condition that $\delta$ satisfies, there holds $2 \left|m_{j}B_{\alpha_{j}}^{-1}\right||z|^{-\alpha_{j}}<\frac{1}{2}$.
Thus,
\begin{displaymath}
\left|u_j\right|\leq 4|z|^{-\alpha_{j}}\left|v_j\right|, ~~j=1,2,
\end{displaymath}
which implies the desired estimates.
\end{pf}
\subsection{The proof of Lemma \ref{lemm:lemma2.6}}

\begin{pf}
For all $z\in\Sigma_{\theta,\delta}$, let
\begin{displaymath}
\left\{\left(\big(z+\rho U(1))^{\alpha_{1}}-m_{1}B_{\alpha_{1}}^{-1}\right)\left((z+\rho U(2))^{\alpha_{2}}-m_{2}B_{\alpha_{2}}^{-1}\right)- m_{1}m_{2}B_{\alpha_{1}}^{-1}B_{\alpha_{2}}^{-1}\right\}u=v.
\end{displaymath}
After simple calculations, we have
\begin{displaymath}\begin{aligned}
u=&\left((z+\rho U(1))^{\alpha_{1}}-m_{1}B_{\alpha_{1}}^{-1}\right)^{-1}\left((z+\rho U(2))^{\alpha_{2}}-m_{2}B_{\alpha_{2}}^{-1}\right)^{-1}v\\
  &+m_{1}m_{2}B_{\alpha_{1}}^{-1}B_{\alpha_{2}}^{-1}\left((z+\rho U(1))^{\alpha_{1}}-m_{1}B_{\alpha_{1}}^{-1}\right)^{-1}
  \left((z+\rho U(2))^{\alpha_{2}}-m_{2}B_{\alpha_{2}}^{-1}\right)^{-1}u.
\end{aligned}
\end{displaymath}
Taking modulus on both sides of the above equality and by Lemma \ref{lemm:lemma2.5}, there is
\begin{displaymath}
|u|\leq 16 |z|^{-\alpha_{1}-\alpha_{2}}|v|+16 \left|m_{1}m_{2}B_{\alpha_{1}}^{-1}B_{\alpha_{2}}^{-1}\right||z|^{-\alpha_{1}-\alpha_{2}}|u|.
\end{displaymath}
According to the condition that $\delta$ satisfies, there holds
$16\left|m_{1}m_{2}B_{\alpha_{1}}^{-1}B_{\alpha_{2}}^{-1}\right||z|^{-\alpha_{1}-\alpha_{2}}<\frac{1}{2}$. Thus,
\begin{displaymath}
|u|\leq 32|z|^{-\alpha_{1}-\alpha_{2}}|v|,
\end{displaymath}
which implies the desired estimate on $H(z)$.

For $H_{\alpha_{1}}(z)$, similarly, let
\begin{displaymath}
\left((z+\rho U(1))^{\alpha_1}-m_1 B_{\alpha_{1}}^{-1}\right)u=v+m_{1}m_{2} B_{\alpha_{1}}^{-1} B_{\alpha_{2}}^{-1}\left((z+\rho U(2))^{\alpha_2}-m_2 B_{\alpha_{2}}^{-1}\right)^{-1}u.
\end{displaymath}
There exists
\begin{displaymath}
u= m_{1}m_{2} B_{\alpha_{1}}^{-1} B_{\alpha_{2}}^{-1}\left((z+\rho U(1))^{\alpha_1}-m_1 B_{\alpha_{1}}^{-1}\right)^{-1}\left((z+\rho U(2))^{\alpha_2}-m_2
   B_{\alpha_{2}}^{-1}\right)^{-1}u
  +\left((z+\rho U(1))^{\alpha_1}-m_1 B_{\alpha_{1}}^{-1}\right)^{-1}v.
\end{displaymath}
Taking modulus on both sides of the above equality and based on Lemma \ref{lemm:lemma2.5}, there holds
\begin{displaymath}
|u|\leq 4 |z|^{-\alpha_{1}}|v|+16 \left|m_{1}m_{2}B_{\alpha_{1}}^{-1}B_{\alpha_{2}}^{-1}\right||z|^{-\alpha_{1}-\alpha_{2}}|u|.
\end{displaymath}
According to the previous estimate, there is $16\left|m_{1}m_{2}B_{\alpha_{1}}^{-1}B_{\alpha_{2}}^{-1}\right||z|^{-\alpha_{1}-\alpha_{2}}<\frac{1}{2}$. Thus,
\begin{displaymath}
|u|\leq 8|z|^{-\alpha_{1}}|v|,
\end{displaymath}
which implies the desired estimate on $H_{\alpha_1}(z)$. By analogy, one can obtain the estimate on $H_{\alpha_2}(z)$.
\end{pf}

\subsection{The proof of Lemma \ref{lemm:lemma2.8}}

\begin{pf}
According to the condition that $\delta$ satisfies and  Lemma \ref{lemm:lemma2.4} and Lemma \ref{lemm:lemma2.6}, the conclusions can be similarly obtained.
\end{pf}

\subsection{The proof of Theorem \ref{the:theorem2.11}}
\begin{pf}
By the inverse Laplace transform, the solution in (\ref{eq:DIVFeynman-kac}) is
\begin{displaymath}\begin{aligned}
 G_{1}(t)=\frac{1}{2\pi i}\int_{\Gamma_{\theta,\delta}}e^{zt}\left(H_{\alpha_{1}}(z)(z+\rho U(1))^{\alpha_{1}-1}G_{1,0}
          -m_{2}B_{\alpha_{2}}^{-1}H(z)(z+\rho U(1))^{\alpha_{1}-1}G_{2,0}\right)dz.
\end{aligned}\end{displaymath}
Taking  $q$-th ($q=0,1$) derivative leads to
\begin{displaymath}\begin{aligned}
 \frac{\partial^{q}}{\partial t^{q}}G_{1}(t)=\frac{1}{2\pi i}\int_{\Gamma_{\theta,\delta}}
           z^{q}e^{zt}\left(H_{\alpha_{1}}(z)(z+\rho U(1))^{\alpha_{1}-1}G_{1,0}
           -m_{2}B_{\alpha_{2}}^{-1}H(z)(z+\rho U(1))^{\alpha_{1}-1}G_{2,0}\right)dz.
\end{aligned}\end{displaymath}
From Lemma \ref{lemm:lemma2.8}, there exists
\begin{displaymath}\begin{aligned}
\left| \frac{\partial^{q}}{\partial t^{q}}G_{1}(t)\right|
   &=\left|\frac{1}{2\pi i}\int_{\Gamma_{\theta,\delta}}z^{q}e^{zt}\left(H_{\alpha_{1}}(z)(z+\rho U(1))^{\alpha_{1}-1}G_{1,0}
              -m_{2}B_{\alpha_{2}}^{-1}H(z)(z+\rho U(1))^{\alpha_{1}-1}G_{2,0}\right)dz\right|\\
   &\leq\frac{8}{\pi}\int_{\Gamma_{\theta,\delta}}e^{Re(z)t}|z|^{q-1}|G_{1,0}||dz|
                +\frac{32\left|m_{2}B_{\alpha_{2}}^{-1}\right|}{\pi}\int_{\Gamma_{\theta,\delta}}e^{Re(z)t}|z|^{q-\alpha_{2}-1}|G_{2,0}||dz|\\
   &\leq\frac{8}{\pi}\left(\int^{+\infty}_{\delta}e^{rt\cos(\theta)}r^{q-1}dr
                +\int^{\theta}_{-\theta}e^{\delta t\cos(\varphi)}\delta^{q}d\varphi\right)|G_{1,0}|\\
   &\quad +\frac{32\left|m_{2}B_{\alpha_{2}}^{-1}\right|}{\pi}\left(\int^{+\infty}_{\delta}e^{rt\cos(\theta)}r^{q-\alpha_{2}-1}dr +
                  \int^{\theta}_{-\theta}e^{\delta t\cos(\varphi)}\delta^{q-\alpha_{2}}d\varphi\right)|G_{2,0}|.
\end{aligned}\end{displaymath}
Let $ rt = s $. Since $t>1/\delta$, there exists
\begin{displaymath}\begin{aligned}
\left| \frac{\partial^{q}}{\partial t^{q}}G_{1}(t)\right|
\leq &\frac{8}{\pi}\left( t^{-q}\int^{+\infty}_{1}e^{s\cos(\theta)}s^{q-1}ds
       +\delta^{q}\int^{\theta}_{-\theta}e^{\delta t\cos(\varphi)}d\varphi\right)|G_{1,0}|\\
     &+\frac{32\left|m_{2}B_{\alpha_{2}}^{-1}\right|}{\pi}\left( t^{\alpha_{2}-q}\int^{+\infty}_{1}e^{s\cos(\theta)}s^{q-\alpha_{2}-1}ds +
       \delta^{q-\alpha_{2}}\int^{\theta}_{-\theta}e^{\delta t\cos(\varphi)}d\varphi\right)|G_{2,0}|\\
\leq &\frac{8}{\pi}\left((-1/\cos(\theta)) t^{-q}e^{\cos(\theta)}+2\theta\delta^{q}e^{\delta T}\right)|G_{1,0}|
\\
&+\frac{32\left|m_{2}B_{\alpha_{2}}^{-1}\right|}{\pi}\left((-1/\cos(\theta))t^{\alpha_{2}-q}e^{\cos(\theta)}+2\theta\delta^{q-\alpha_{2}}e^{\delta T}\right)|G_{2,0}|.
\end{aligned}\end{displaymath}
Similarly,
\begin{displaymath}\begin{aligned}
\left|G_{2}^{(q)}(t)\right|
          \leq \frac{32\left|m_{1}B_{\alpha_1}^{-1}\right|}{\pi}\left((-1/\cos(\theta))t^{\alpha_1-q}e^{\cos(\theta)}+2\theta\delta^{q-\alpha_1}e^{\delta T}\right)|G_{1,0}|
 +\frac{8}{\pi}\left((-1/\cos(\theta))t^{-q}e^{\cos(\theta)}+2\theta\delta^{q}e^{\delta T}\right)|G_{2,0}|.\\
\end{aligned}\end{displaymath}
The proof is completed.
\end{pf}

\section{The proofs of the results presented in Section \ref{subsec:NCIM}}

\subsection{\textbf{The proof of Proposition \ref{Thm:Prop1}}}
\label{det:strip}
\begin{pf}
Let the denominator of $H(z)$ be zero, i.e.,
\begin{equation}\label{eq:sigular}
\left((z+\rho U(1))^{\alpha_{1}}-C_{1}\right)\left((z+\rho U(2))^{\alpha_{2}}-C_{2}\right)-C_{1}C_{2}=0.
\end{equation}
Denote $u_{1}+iv_{1}:=(z+\rho U(1))^{\alpha_{1}}$, $u_{2}+iv_{2}:=(z+\rho U(2))^{\alpha_{2}}$. Then (\ref{eq:sigular}) can be rewritten as
\begin{equation}\label{eq:esigular}
\left((u_{1}-C_{1})(u_{2}-C_2)-v_{1}v_{2}-C_{1}C_{2}\right)+i\left((u_{1}-C_{1})v_{2}+(u_{2}-C_{2})v_{1}\right)=0.
\end{equation}
In the sequence, we prove that if (\ref{condition01}) holds, (\ref{eq:esigular}) does not hold.

We divide the proof into the following cases:
\begin{description}
  \item[Case I:] For the case $u_1>2 C_1$ and $u_2>2 C_2$, if $v_{1}v_{2} \leq 0$, then $(u_{1}-C_{1})(u_{2}-C_2)-v_1v_2-C_{1}C_{2}\neq0$; otherwise, $(u_{1}-C_{1})v_{2}+(u_{2}-C_{2})v_{1} \neq 0$.
  \item[Case II:] For the case $|v_1|>|C_1|$ and $|v_2|>|C_2|$,
  \begin{description}
    \item[$(a)$] when $u_1\leq2 C_1$ and $u_2\leq2 C_2$:
        If $(u_{1}-C_{1})(u_{2}-C_2)>0$, for $v_1  v_2 >0$, there is $(u_{1}-C_{1})v_{2}+(u_{2}-C_{2})v_{1} \neq 0$;
           for $v_1  v_2 < 0$, we have $(u_{1}-C_{1})(u_{2}-C_2)-v_1v_2-C_{1}C_{2}\neq0$.
        If $(u_{1}-C_{1})(u_{2}-C_2)=0$, $(u_{1}-C_{1})v_{2}+(u_{2}-C_{2})v_{1} \neq 0$.
        If $(u_{1}-C_{1})(u_{2}-C_2)<0$, for $v_1  v_2 >0$, there is $(u_{1}-C_{1})(u_{2}-C_2)-v_1v_2-C_{1}C_{2}\neq0$;
           for $v_1  v_2 < 0$, we have $(u_{1}-C_{1})v_{2}+(u_{2}-C_{2})v_{1} \neq 0$.
    \item[$(b)$] when $(u_1-2 C_1)(u_2-2 C_2)\leq0$:
        If $(u_1-C_1)(u_2-C_2)>0$, for $v_1v_2>0$, there is  $(u_{1}-C_{1})v_{2}+(u_{2}-C_{2})v_{1} \neq 0$;
           for $v_1v_2<0$, we have  $(u_{1}-C_{1})(u_{2}-C_2)-v_1v_2-C_{1}C_{2}\neq0$.
        If $(u_1-C_1)(u_2-C_2)=0$, there is  $(u_{1}-C_{1})(u_{2}-C_2)-v_1v_2-C_{1}C_{2}\neq0$.
        If $(u_{1}-C_{1})(u_{2}-C_2)<0$, for $v_1  v_2 >0$, there is $(u_{1}-C_{1})(u_{2}-C_2)-v_1v_2-C_{1}C_{2}\neq0$.
  \end{description}
\end{description}
To sum up, when (\ref{condition01}) holds, $H(z)$ is analytic.
\end{pf}

\subsection{\textbf{The proof of (\ref{condition02})}}
\label{condi}
\begin{pf}
For the proof of (\ref{condition02}), we spilt it into the following two cases:
\begin{description}
  \item[Case I:] Let $\textmd{Re}(z+\rho U(j))>0$,  i.e., $\textmd{Re}(z)>-\textmd{Re}(\rho U(j))$. For fixed $0<\alpha_j<1$, there are $\theta_j:=\arg(z+\rho U(j))\in(-\frac{\pi}{2}, \frac{\pi}{2})$ and $\alpha_j\theta_j\in(-\frac{\alpha_j\pi}{2},\frac{\alpha_j\pi}{2})$. Further, we have $\textmd{Re}\left((z+\rho U(j)\right)^{\alpha_j})=|z+\rho U(j)|^{\alpha_j}\cos(\alpha_j\theta_j)>|z+\rho U(j)|^{\alpha_j}\cos(\frac{\alpha_j\pi}{2})$. If $|z+\rho U(j)|^{\alpha_j}\cos(\frac{\alpha_j\pi}{2})>2|C_j|$, that is $|z+\rho U(j)|>(2|C_j|/\cos(\alpha_j\pi/2))^{1/\alpha_j}$, then $\textmd{Re}(z+\rho U(j))=|z+\rho U(j)|\cos(\alpha_j\theta_j)>2|C_j|/\cos(\alpha_j\pi/2))^{1/\alpha_j}\cos(\alpha_j\pi)$. Hence, denote
      \begin{displaymath}
        d_1^{\langle j\rangle}:=\left(2|C_j|/\cos(\alpha_j\pi/2)\right)^{1/\alpha_j}-\textmd{Re}(\rho U(j)),~ ~ j=1,2.
      \end{displaymath}
      Once $\textmd{Re}(z)>d_1:=\max\{d_1^{\langle 1\rangle}, d_1^{\langle 2\rangle}\}$, there holds $\textmd{Re}\left((z+\rho U(j)\right)^{\alpha_j})>2|C_j|$, $j=1,2$. By (\ref{condition01}), $H(z)$ is analytic.

  \item[Case II:] By (\ref{condition01}), $H(z)$ is analytic if $\left|\textmd{Im}((z+\rho U(j))^{\alpha_j})\right|>|C_j|$ and $|\textmd{Re}\left((z+\rho U(j)\right)^{\alpha_j})|\leq2|C_j|$ at the same time. So, if $\textmd{Im}((z+\rho U(j))^{\alpha_j})=|z+\rho U(i)|^{\alpha_j}\sin(\alpha_j\theta_j)>|z+\rho U(j)|^{\alpha_j}\sin((|C_j|^{2}/(5|C_j|^{2}))^{1/2})>|C_j|$ (where $\theta_j$ are defined as above), i.e., $|z+\rho U(i)|>(\sqrt{5}|C_j|)^{1/\alpha_j}$, then $H(z)$ is analytic. With this, let $|z+\rho U(j)|>|\textmd{Im}(z+\rho U(j))|\geq\left|~|\textmd{Im}(z)|-|\textmd{Im}(\rho U(j))|~\right|>(\sqrt{5}|C_j|)^{1/\alpha_j}$, and denote
      \begin{displaymath}
        d_2^{\langle j\rangle}:=\left(\sqrt{5}|C_j|)\right)^{1/\alpha_j}+|\textmd{Im}(\rho U(j))|,~ ~ j=1,2.
      \end{displaymath}
      Then, for $|\textmd{Im}(z)|>d_2:=\max\{d_2^{\langle 1\rangle}, d_2^{\langle 2\rangle}\}$, there is $\left|\textmd{Im}((z+\rho U(j))^{\alpha_j})\right|>|C_j|$, and $H(z)$ is analytic.
\end{description}
Above all, if $z$ satisfies the conditions in  (\ref{condition02}), then $H(z)$ is analytic.
\end{pf}

\subsection{\textbf{Determination of the open strips $S^{\langle1\rangle}$ and $S^{\langle2\rangle}$}}
\label{det:strips}
\begin{figure}[t!]
  \centering
  \includegraphics[width=0.32\textwidth]{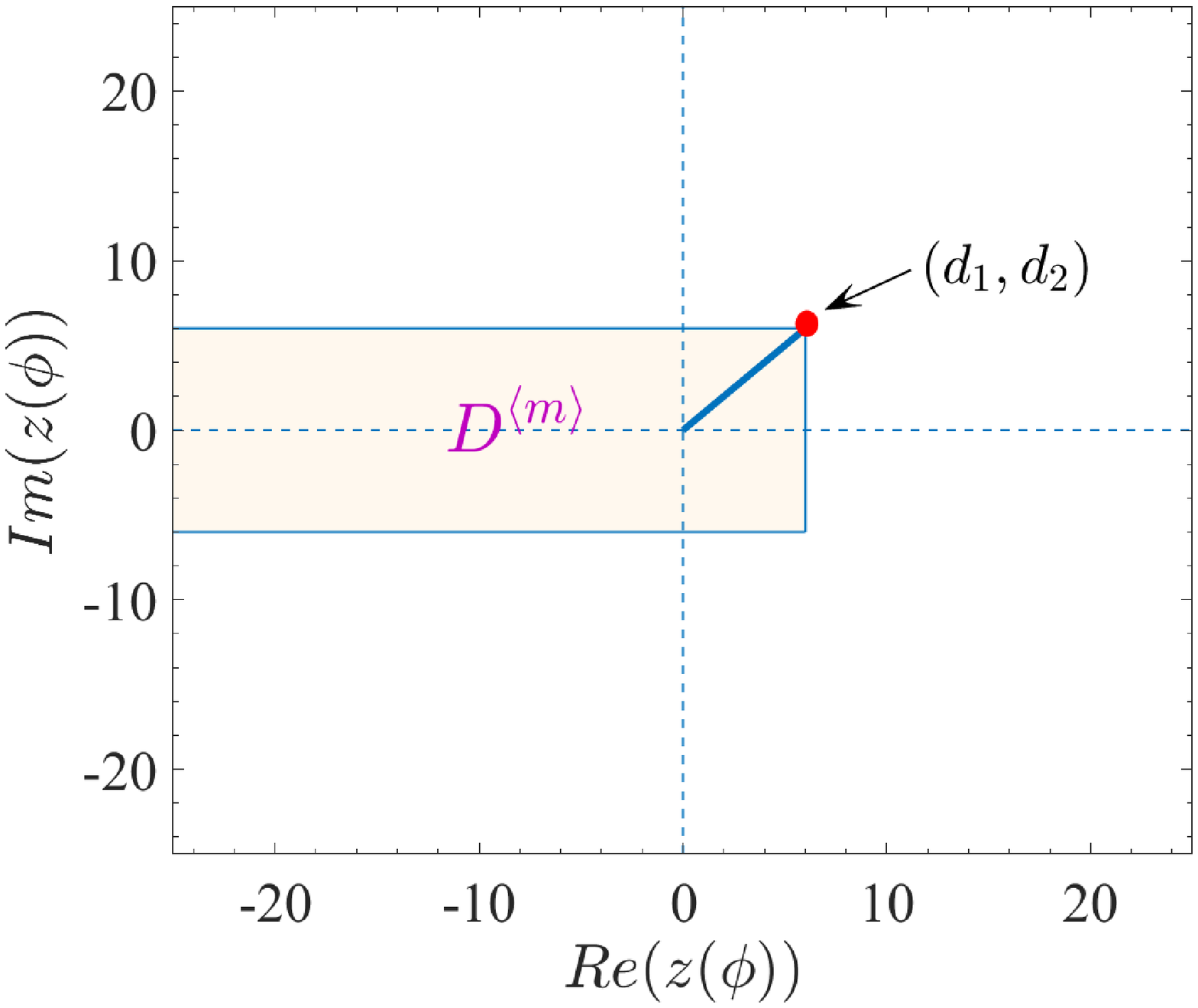}
  \includegraphics[width=0.32\textwidth]{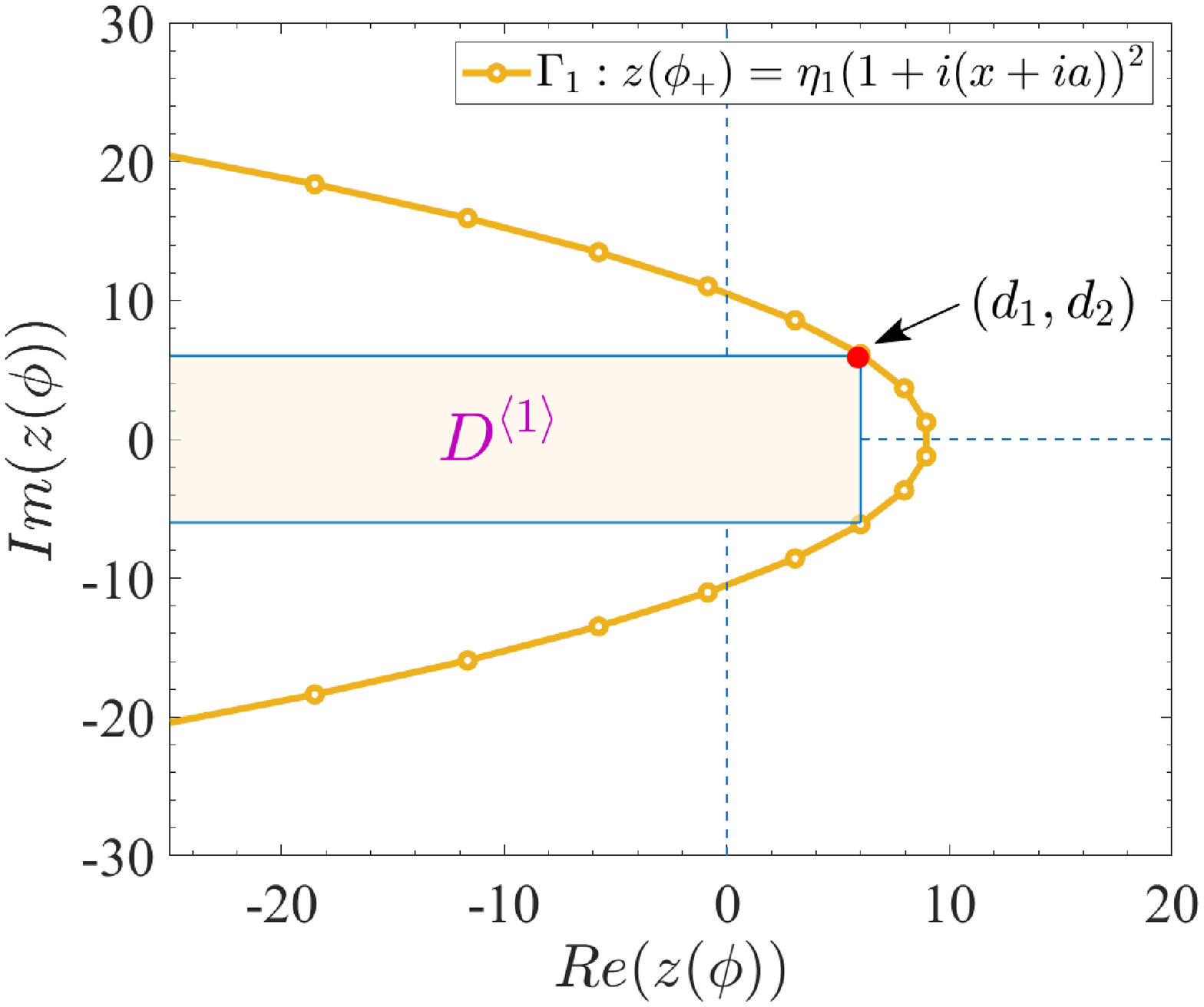}
  \includegraphics[width=0.32\textwidth]{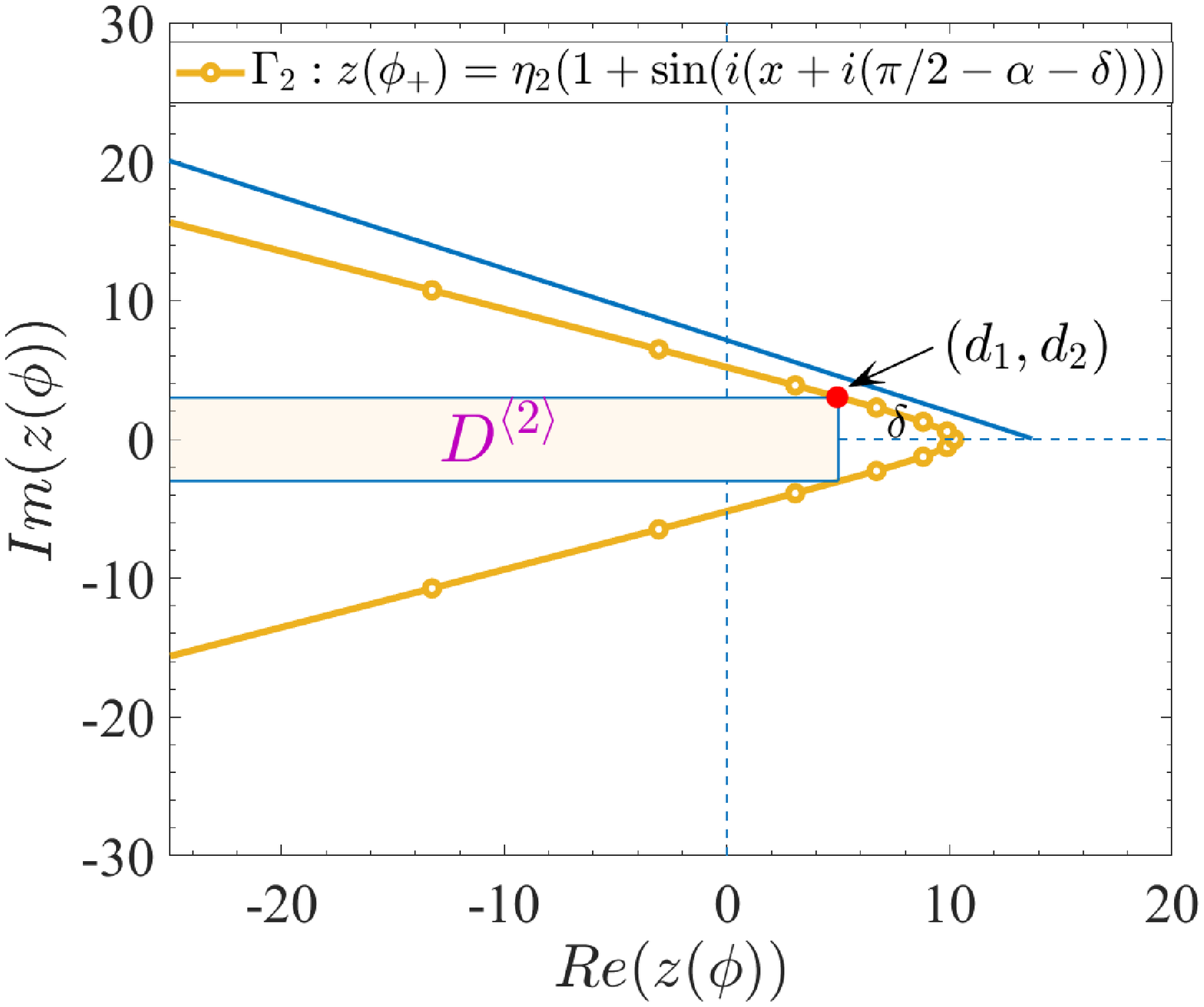}
\caption{ \emph{Diagrammatic sketch for the condition (\ref{condition02}) and the left boundary of $N^{\langle m\rangle}$ for different integral contours}. (\emph{left}) is a domain determined by (\ref{condition02}); (\emph{center}) is to determine the parameters of the open strip $S^{\langle1\rangle}$ when the CIM is with the parabolic contour  $\Gamma_{1}$; (\emph{right}) is to determine the parameters of the open strip $S^{\langle2\rangle}$  when the CIM is with the hyperbolic contour $\Gamma_{2}$.}
\label{Fig:conditons}
\end{figure}

For strip $S^{\langle1\rangle}$, which corresponds to the integral contour $\Gamma_1$, as shown in Figure \ref{Fig:conditons} (\emph{center}). Denote $\phi_+:=x+ia\in S$ with $a>0$. Then
\begin{equation}\label{eq:pa}
z(\phi_+)=\eta_1\left((1-a)^{2}-x^{2}\right)+2i\eta_1x(1-a),
\end{equation}
$z(\phi_+)=z(x+ia)$ is the left boundary of $N^{\langle m\rangle}$. which can be expressed as the parabola
\begin{equation}\label{eq:pa1}
u=\eta_1\left((1-a)^{2}-\frac{v^{2}}{4\eta_1^{2}(1-a)^{2}}\right),
\end{equation}
if denoting $z=u+iv$. As $a$ increases from $0$ to $1$, the parabola (\ref{eq:pa}) closes and reduces to the negative real axis. Hence, the left boundary of $N_r$ determines the maximum value of $a$, which means that the parabola (\ref{eq:pa1}) passing through the point $(d_{1},d_{2})$, yields
\begin{equation}\label{dete:a}
    a=1-\sqrt{\frac{d_{1}+\sqrt{d_{1}^{2}+d_{2}^{2}}}{2\eta_{1}}},
\end{equation}
where $d_1$ and $d_2$ are given in \ref{condi}.

Denote $\phi_{-}:=x-id$ with $d>0$. The image of this horizontal line is
\begin{equation}\label{eq:pa2}
z(\phi_{-})=\eta_1\left((1+d)^{2}-x^{2}\right)+2i\eta_{1}x(1+d).
\end{equation}
As $d$ away from $0$, the parabola (\ref{eq:pa2}) widens and moves to the right.  The optimal $d$ and $\eta_1$  determined in (\ref{eq:d}) and (\ref{eq:eta_1}), respectively.

For the strip $S^{\langle2\rangle}$, which corresponds to the integral contour $\Gamma_2$, as shown in Figure \ref{Fig:conditons} (\emph{right}). From the expression of $\Gamma_2$, the image of the horizontal line $\phi=x+iy$ is
\begin{equation}\label{eq:hz}
z(\phi)=\eta_2\left(1-\sin(\alpha+y)\cosh(x))+i\eta_2\cos(\alpha+y)\sinh(x)\right),
\end{equation}
which can be expressed as the hyperbola
\begin{equation}
\left(\frac{\eta_2-u}{\sin(\alpha+y)}\right)^{2}-\left(\frac{v}{\cos(\alpha+y)}\right)^{2}=\eta_{2}^{2},
\end{equation}
if denoting $z=u+iv$. As $r$ increases from $0$ to $r=\pi/2-\alpha$, the left branch of parabola (\ref{eq:hz}) closes and degenerates into the negative real axis. While, for $y<0$, when $r$ decreases from $0$ to $-\alpha$, the hyperbola widens and becomes a vertical line. The minimum value of $\delta$ can be obtained by taking $z(\phi_+):=z(x+i(\pi/2-\alpha-\delta))$ as the left boundary of $N^{\langle m\rangle}$. More specifically, $z(\phi_+)$ can be expressed as a hyperbola with the asymptotes
\begin{equation}
v=\pm\tan(\delta)(\eta_2-u),
\end{equation}
if denoting $z=u+iv$.  Let one of the asymptotic lines $v=\tan(\delta)(\eta_2-u)$ pass through the fixed point $(d_1,d_2)$, which results in  $\delta=\arctan\left(\frac{d_{2}}{\eta_{2}-d_{1}}\right)$. Besides, the optimal parameter $\alpha$ and $\eta_2$ are determined by maximizing $Q(\alpha)=\frac{\pi^{2}-2\pi\alpha-2\pi\delta}{A(\alpha)}$ (see Subsection \ref{sec:opt} and (\ref{eq:eta_2})). Then, the strip $S^{\langle2\rangle}$ is determined.

\subsection{\textbf{The proof of Lemma \ref{lemm:lemma3.1}}}
\label{proof:lemm3}
\begin{pf}
For the first part of the lemma, as $z(\phi)=\eta_{1}(i\phi+1)^{2}$ and $z'(\phi)=i2\eta_{1}(i\phi+1)$, from the expression of $v_{1}(t,\phi)$, there is
\begin{displaymath}
v_{1}(t,\phi)=\frac{\eta_{1}}{\pi}e^{\eta_{1}(i\phi+1)^{2}t}(1+i\phi)\widehat{G}_{1}\left(\eta_{1}(1+i\phi)^{2}\right)~ ~ \forall~ \phi\in S^{\langle1\rangle}.
\end{displaymath}

Choose $\phi_+=x+iy\in S^{\langle1\rangle}$, $0<y<a<1$, from the upper half plane of the strip $S^{\langle1\rangle}$, for $t>0$, there holds
\begin{displaymath}
v_{1}(t,x+iy)
=\frac{\eta_{1}e^{\eta_{1}\left((1-y)^{2}-x^{2}\right)t}}{\pi}e^{i2\eta_{1}tx(1-y)}(1-y+ix)\widehat{G}_{1}\left(\eta_{1}(1-y+ix)^{2}\right).
\end{displaymath}
Taking the modulus of the left and right sides of the above formula leads to
\begin{displaymath}
 |v_{1}(t,x+iy)|
\leq \frac{\eta_{1}e^{\eta_{1}(1-y)^{2}t}}{\pi}e^{-x^{2}\eta_{1}t}\left|(1-y+ix)\widehat{G}_{1}\left(\eta_{1}(1-y+ix)^{2}\right)\right|.
\end{displaymath}
Denote $l=1-y$, $l\in(1-a,1)$. From Lemma \ref{lemm:lemma2.8}, there holds $\left|\widehat{G}_{1}(z)\right|\leq 16|z|^{-1}|G_{1,0}|+64\left|m_1B_{\alpha_1}^{-1}\right||z|^{-1-\alpha_{2}}|G_{2,0}|$. With these, we have
\begin{small}\begin{displaymath}
     \left|(l+ix)\widehat{G}_{1}\left(\eta_{1}(l+ix)^{2}\right)\right|
\leq  \left(\frac{16|l+ix|}{|\eta_{1}(l+ix)^{2}|}|G_{1,0}|
+\frac{64\left|m_1B_{\alpha_1}^{-1}\right||l+ix|}{|\eta_{1}(l+ix)^{2}|^{1+\alpha_{2}}}|G_{2,0}|\right).
\end{displaymath}\end{small}
Since
\begin{displaymath}
\frac{|l+ix|}{|\eta_{1}(l+ix)^{2}|} \leq \frac{1}{\eta_{1}(1-a)}~{\rm and}~
\frac{|l+ix|}{|\eta_{1}(l+ix)^{2}|^{1+\alpha_{2}}}\leq \frac{1}{\eta_{1}^{1+\alpha_{2}}(1-a)^{1+2\alpha_{2}}},
\end{displaymath}
we have
\begin{displaymath}\begin{aligned}
|v_{1}(t,x+ir)|
\leq &\frac{\eta_{1}e^{\eta_{1}t}}{\pi}\left(\frac{16}{\eta_{1}(1-a)}|G_{1,0}|
                       +\frac{64\left|m_1B_{\alpha_1}^{-1}\right|}{\eta_{1}^{1+\alpha_{2}}(1-a)^{1+2\alpha_{2}}}|G_{2,0}|\right)e^{-\eta_{1}tx^{2}}\\
\leq &\frac{e^{\eta_{1}t}}{\pi}\mathrm{max}\left\{16(1-a)^{-1},64\left|m_1B_{\alpha_1}^{-1}\right|(1-a)^{-(1+2\alpha_2)}\right\}
                       \left(|G_{1,0}|+\eta_{1}^{-\alpha_{2}}|G_{2,0}|\right)e^{-\eta_{1}tx^{2}}.
\end{aligned}\end{displaymath}

Choose $\phi_-=x-iy\in S^{\langle1\rangle}$, $-d<-y$, from the lower half plane of $S^{\langle1\rangle}$, for $t>0$, there holds
\begin{displaymath}\begin{aligned}
|v_{1}(t,x-iy)|
&\leq \frac{\eta_1e^{\eta_{1}(1+d)^{2}t}}{\pi}\left(16\eta_{1}^{-1}|G_{1,0}|
     +64\left|m_1B_{\alpha_1}^{-1}\right|\eta_{1}^{-(1+\alpha_{2})}|G_{2,0}|\right)e^{-\eta_{1}tx^{2}}\\
&\leq \frac{e^{\eta_{1}t}}{\pi}\mathrm{max}\left\{16e^{\left(2d+d^2\right)\eta_1t},
          64\left|m_1B_{\alpha_1}^{-1}\right|e^{\left(2d+d^2\right)\eta_1t}\right\}
\left(|G_{1,0}|+\eta_{1}^{-\alpha_{2}}|G_{2,0}|\right)e^{-\eta_{1}tx^{2}}.
\end{aligned}\end{displaymath}
Therefore, for $\phi\in S^{\langle1\rangle}$, by denoting
\begin{equation}\label{con:B}
 B:=\mathrm{max}\left\{16(1-a)^{-1},
      64\left|m_1B_{\alpha_1}^{-1}\right|(1-a)^{-(1+2\alpha_2)},
      16e^{\left(2d+d^2\right)\eta_{1}t},
      64\left|m_1B_{\alpha_1}^{-1}\right|e^{\left(2d+d^2\right)\eta_{1}t}
 \right\},
\end{equation}
we have
\begin{displaymath}
|v_{1}(t,\phi)| \leq \frac{Be^{\eta_{1}t}}{\pi}\left(|G_{1,0}|
+\eta_{1}^{-\alpha_{2}}|G_{2,0}|\right)e^{-\eta_{1}tx^{2}}.
\end{displaymath}

For the case of $S^{\langle2\rangle}$, the proof is similar to the previous one. Hence $v_{1}(t,\phi)$ is analytical in the strip $S^{\langle2\rangle}$ and $z'(\phi)=i\eta_2\cos(i\phi-\alpha)$. By choosing $\phi_{+}=x+iy\in S^{\langle2\rangle}$, $0<y<\pi/2-\delta-\alpha$, form the upper half plane of the strip $S^{\langle2\rangle}$, for
$t>0$, there holds
\begin{displaymath}
v_{1}(t,x+iy)
=\frac{\eta_{2}}{2\pi}e^{\eta_{2}(1-\sin(\alpha+y-ix))t}\cos(\alpha+y-ix)\widehat{G}_{1}(\eta_{2}(1-\sin(\alpha+y-ix))).
\end{displaymath}
Denote $l'=\alpha+y$, $l'\in(\alpha, \pi/2-\delta)$. After taking the modulus of the left and right sides of the above formula, we get
\begin{displaymath}
|v_{1}(t,x+ir)|
\leq \frac{\eta_{2}}{2\pi}e^{\eta_{2}t(1-\sin l'\cosh x)}\left|\cos(l'-ix)\widehat{G}_{1}(\eta_{2}(1-\sin(l'-ix)))\right|.
\end{displaymath}
Moreover, since
\begin{displaymath}\begin{aligned}
     \left|\cos(l'-ix)\widehat{G}_{1}\left(\eta_{2}(1-\sin (l'-ix))\right)\right|
\leq  \frac{16|\cos(l'-ix)|}{|\eta_{2}(1-\sin(l'-ix))|}|G_{1,0}|
   +\frac{64\left|m_1B_{\alpha_1}^{-1}\right||\cos(l'-ix)|}{|(\eta_{2}(1-\sin(l'-ix)))^{1+\alpha_{2}}|}|G_{2,0}|,
\end{aligned}\end{displaymath}
\begin{displaymath}
\frac{|\cos(l'-ix)|}{|\eta_{2}(1-\sin(l'-ix))|}\leq \frac{1}{\eta_{2}}\sqrt{\frac{1+\sin l'}{1-\sin l'}},
\end{displaymath}
and
\begin{displaymath}\begin{aligned}
\frac{|\cos(l'-ix)|}{|(\eta_{2}(1-\sin(l'-ix)))^{1+\alpha_{2}}|}
 &=\frac{1}{|(\eta_{2}(1-\sin(l'-ix)))^{\alpha_{2}}|}\frac{|\cos(l'-ix)|}{|\eta_{2}(1-\sin(l'-ix))|}\\
 &\leq \frac{1}{\eta_{2}^{1+\alpha_{2}}(\cosh(x)-\sin{l'})^{\alpha_{2}}}\sqrt{\frac{1+\sin l'}{1-\sin l'}},
\end{aligned}\end{displaymath}
then, by $\sin\alpha<\sin l'=\sin(\alpha+y)$, there holds
\begin{displaymath}\begin{aligned}
\left|v_{1}(t,x+ir)\right|
\leq & \frac{8e^{\eta_{2}t}}{\pi}\sqrt{\frac{1+\sin(\alpha+y)}{1-\sin(\alpha+y)}}\left(|G_{1,0}|
     +\frac{4\left|m_1B_{\alpha_1}^{-1}\right|}{\eta_{2}^{\alpha_{2}}(1-\sin(\alpha+y))^{\alpha_{2}}}|G_{2,0}|\right)
      e^{-\eta_{2}t\sin{\alpha}\cosh x}.
\end{aligned}\end{displaymath}

Choosing $\phi_-=x-iy\in S^{\langle2\rangle}$, $-\alpha<-y$, from the lower half plane of $S^{\langle2\rangle}$, for $t>0$, we deduce
\begin{displaymath}\begin{aligned}
\left|v_{1}(t,x-iy)\right|
\leq & \frac{8e^{\eta_{2}t}}{\pi}\sqrt{\frac{1+\sin(\alpha+y)}{1-\sin(\alpha+y)}}\left(|G_{1,0}|
     +\frac{4\left|m_1B_{\alpha_1}^{-1}\right|}{\eta_{2}^{\alpha_{2}}(1-\sin(\alpha+y))^{\alpha_{2}}}|G_{2,0}|\right)
      e^{-\eta_{2}t\sin{\alpha}\cosh x}.
\end{aligned}\end{displaymath}
Therefore, for $\phi\in S^{\langle2\rangle}$, $0<y<\min\{\alpha, \pi/2-\alpha-\delta\}$, there holds
\begin{displaymath}\begin{aligned}
\left|v_{1}(t,\phi)\right|
\leq & \frac{8e^{\eta_{2}t}}{\pi}\sqrt{\frac{1+\sin(\alpha+y)}{1-\sin(\alpha+y)}}\left(|G_{1,0}|
     +\frac{4\left|m_1B_{\alpha_1}^{-1}\right|}{\eta_{2}^{\alpha_{2}}(1-\sin(\alpha+y))^{\alpha_{2}}}|G_{2,0}|\right)
      e^{-\eta_{2}t\sin{\alpha}\cosh x}.
\end{aligned}\end{displaymath}

Similar estimates for $v_{2}(t,\phi)$ can be obtained on the strips $S^{\langle1\rangle}$ and $S^{\langle2\rangle}$, respectively.
\end{pf}

\subsection{\textbf{The proof of Theorem \ref{thm:estev}}}
\label{proof:thm}
\begin{pf}
We firstly choose the integral contour $z(\phi)$ as defined in (\ref{eq:paraboliccontour}) with the uniform step-size $h_1$.
For $t>0$, by (\ref{est:bound1}), there holds
\begin{displaymath}
\left|G^{\langle 1\rangle}_{1,N}(t)\right| \leq h_{1}\sum\limits_{|k|\leq N-1}|v_{1}(t,kh_{1})|
\leq \frac{2 \eta_{1}Be^{\eta_{1} t}}{\pi}\left(\eta_{1}^{-1}|G_{1,0}|
+\eta_{1}^{-(1+\alpha_{2})}|G_{2,0}|\right)\int^{(N-1)h_{1}}_{0}e^{-x^{2}\eta_{1}t}dx.
\end{displaymath}
Since
\begin{displaymath}
 \int^{(N-1)h_{1}}_{0}e^{-x^{2}\eta_{1}t}dx\leq\int^{\infty}_{0}e^{-x^{2}\eta_{1}t}dx=\frac{\sqrt{\eta_{1}t\pi}}{2\eta_{1}t},
\end{displaymath}
the result of the first part of the theorem holds.

For the rest part of this theorem, we choose the hyperbolic integral contour $z(\phi)$ defined in (\ref{eq:hyperboliccontour}) with the uniform step-size  $h_2$. For $t>0$ and $h_2>0$, by (\ref{est:bound2}), there holds
\begin{displaymath}\begin{aligned}
\left|G^{\langle 2\rangle}_{1,N}(t)\right|
& \leq h_{2}\sum\limits_{|k|\leq N-1}|v_{1}(t,kh_{2})|\\
& \leq \frac{16e^{\eta_{2}t}}{\pi}\sqrt{\frac{1+\sin (\alpha+r)}{1-\sin (\alpha+r)}}\left(|G_{1,0}|
     +\frac{4\left|m_1B_{\alpha_1}^{-1}\right|}{\eta_{2}^{\alpha_{2}}(1-\sin (\alpha+r))^{\alpha_{2}}}|G_{2,0}|\right)\int^{(N-1)h_{2}}_{0}e^{-\eta_{2}t\sin{\alpha}\cosh x}dx.
\end{aligned}\end{displaymath}
Thanks to Lemma \ref{et:lemma}, there holds
\begin{displaymath}
\int^{(N-1)h_{2}}_{0}e^{-\eta_{2}\sin(\alpha)\cosh(x)t}dx\leq \int^{\infty}_{0}e^{-\eta_{2}\sin(\alpha)\cosh(x)t}dx\leq L(\eta_{2}t\sin(\alpha)).
\end{displaymath}

Similarly, we can obtain the stability results about $G^{\langle m\rangle}_{2,N}(t)$, $m=1,2$.
\end{pf}

\section{The time-marching schemes for the Feynman-Kac system (Subsection \ref{subsec:PC})}
\label{sec:Pre-Corre}
Here we provide the time-marching schemes for (\ref{bianxing}) directly. To make things clear for readers, here we demonstrate them in a tedious way. After integrating from $0$ to $t$ on the left and right sides of (\ref{bianxing}), according to (\ref{eq:integral}), for $t=t_{n+1}, n=0, 1, ..., M $, there hold
\begin{align}\label{sys:integral}
G_{1}(t_{n+1})=&\frac{m_{1}B_{\alpha_{1}}^{-1}}{\Gamma(\alpha_{1})}\left(\int_{0}^{t_{n+1}}\frac{e^{-\rho U(1)(t_{n+1}-\tau)}}{(t_{n+1}-\tau)^{1-\alpha_{1}}}
                 G_{1}(\tau)d\tau+(\rho U(1))^{1-\alpha_{1}}\gamma(\alpha_{1},1)\int_{0}^{t_{n+1}}G_{1}(\tau)d\tau\right) \nonumber \\
               -&\frac{m_{2}B_{\alpha_{2}}^{-1}}{\Gamma(\alpha_{2})}\left(\int_{0}^{t_{n+1}}\frac{e^{-\rho U(2)(t_{n+1}-\tau)}}{(t_{n+1}-\tau)^{1-\alpha_{2}}}
                 G_{2}(\tau)d\tau+(\rho U(2))^{1-\alpha_{2}}\gamma(\alpha_{2},1)\int_{0}^{t_{n+1}}G_{2}(\tau)d\tau\right) \nonumber \\
               -&\rho U(1)\int_{0}^{t_{n+1}}G_{1}(\tau)d\tau+G_{1,0},
\end{align}
\begin{align}\label{sys:integral2}
G_{2}(t_{n+1})=&\frac{m_{2}B_{\alpha_{2}}^{-1}}{\Gamma(\alpha_{2})}\left(\int_{0}^{t_{n+1}}\frac{e^{-\rho U(2)(t_{n+1}-\tau)}}{(t_{n+1}-\tau)^{1-\alpha_{2}}}
                 G_{2}(\tau)d\tau+(\rho U(2))^{1-\alpha_{2}}\gamma(\alpha_{2},1)\int_{0}^{t_{n+1}}G_{2}(\tau)d\tau\right) \nonumber \\
              -&\frac{m_{1}B_{\alpha_{1}}^{-1}}{\Gamma(\alpha_{1})}\left(\int_{0}^{t_{n+1}}\frac{e^{-\rho U(1)(t_{n+1}-\tau)}}{(t_{n+1}-\tau)^{1-\alpha_{1}}}
                 G_{1}(\tau)d\tau+(\rho U(1))^{1-\alpha_{1}}\gamma(\alpha_{1},1)\int_{0}^{t_{n+1}}G_{1}(\tau)d\tau\right)  \nonumber \\
               -&\rho U(2)\int_{0}^{t_{n+1}}G_{2}(\tau)d\tau+G_{2,0},
\end{align}
where $\gamma(\alpha_{j},1)$, $j=1,2$ are the incomplete Gamma function.  With these, by the technics mentioned in \cite{Freed02}, one can obtain the following numerical schemes.

\textbf{Case I:} For $n=0$,
\begin{equation}\label{eq:n0C1}
 \begin{aligned}
G_{1}(t_{1})=&
      \frac{m_{1}B_{\alpha_{1}}^{-1}}{\Gamma(\alpha_{1})}\left(\frac{h^{\alpha_{1}}}{\alpha_{1}(\alpha_{1}+1)}\alpha_{1}e^{-\rho U(1)h}
        +\frac{h}{2}(\rho U(1))^{1-\alpha_{1}}\gamma(\alpha_{1},1)\right)G_{1,0}\\
      &+\frac{m_{1}B_{\alpha_{1}}^{-1}}{\Gamma(\alpha_{1})}\left(\frac{h^{\alpha_{1}}}{\alpha_{1}(\alpha_{1}+1)}
        +\frac{h}{2}(\rho U(1))^{1-\alpha_{1}}\gamma(\alpha_{1},1)\right)G_{1}(t_{1})\\
      &-\frac{m_{2}B_{\alpha_{2}}^{-1}}{\Gamma(\alpha_{2})}\left(\frac{h^{\alpha_{2}}}{\alpha_{2}(\alpha_{2}+1)}\alpha_{2}e^{-\rho U(2)h}
        +\frac{h}{2}(\rho U(2))^{1-\alpha_{2}}\gamma(\alpha_{2},1)\right)G_{2,0}\\
      &-\frac{m_{2}B_{\alpha_{2}}^{-1}}{\Gamma(\alpha_{2})}\left(\frac{h^{\alpha_{2}}}{\alpha_{2}(\alpha_{2}+1)}
        +\frac{h}{2}(\rho U(2))^{1-\alpha_{2}}\gamma(\alpha_{2},1)\right)G_{2}(t_{1})\\
      &-\frac{\rho U(1)h-2}{2}G_{1,0}-\frac{\rho U(1)h}{2}G_{1}(t_{1}),
\end{aligned}
\end{equation}
\begin{equation}\label{eq:n0c2}
 \begin{aligned}
G_{2}(t_{1})=
      &\frac{m_{2}B_{\alpha_{2}}^{-1}}{\Gamma(\alpha_{2})}\left(\frac{h^{\alpha_{2}}}{\alpha_{2}(\alpha_{2}+1)}\alpha_{2}e^{-\rho U(2)h}
        +\frac{h}{2}(\rho U(2))^{1-\alpha_{2}}\gamma(\alpha_{2},1)\right)G_{2,0}\\
      &+\frac{m_{2}B_{\alpha_{2}}^{-1}}{\Gamma(\alpha_{2})}\left(\frac{h^{\alpha_{2}}}{\alpha_{2}(\alpha_{2}+1)}
        +\frac{h}{2}(\rho U(2))^{1-\alpha_{2}}\gamma(\alpha_{2},1)\right)G_{2}(t_{1})\\
      &-\frac{m_{1}B_{\alpha_{1}}^{-1}}{\Gamma(\alpha_{1})}\left(\frac{h^{\alpha_{1}}}{\alpha_{1}(\alpha_{1}+1)}\alpha_{1}e^{-\rho U(1)h}
        +\frac{h}{2}(\rho U(1))^{1-\alpha_{1}}\gamma(\alpha_{1},1)\right)G_{1,0}\\
      &-\frac{m_{1}B_{\alpha_{1}}^{-1}}{\Gamma(\alpha_{1})}\left(\frac{h^{\alpha_{1}}}{\alpha_{1}(\alpha_{1}+1)}
        +\frac{h}{2}(\rho U(1))^{1-\alpha_{1}}\gamma(\alpha_{1},1)\right)G_{1}(t_{1})\\
      &-\frac{\rho U(2)h-2}{2}G_{2,0}-\frac{\rho U(2)h}{2}G_{2}(t_{1}).
\end{aligned}
\end{equation}
Denote
\begin{align*}
a_{11}&=\frac{m_{1}B_{\alpha_{1}}^{-1}}{\Gamma(\alpha_{1})}\left(\frac{h^{\alpha_{1}}}{\alpha_{1}(\alpha_{1}+1)}
        +\frac{h}{2}(\rho U(1))^{1-\alpha_{1}}\gamma(\alpha_{1},1)\right)-\frac{\rho U(1)h}{2},\\
a_{12}&=-\frac{m_{2}B_{\alpha_{2}}^{-1}}{\Gamma(\alpha_{2})}\left(\frac{h^{\alpha_{2}}}{\alpha_{2}(\alpha_{2}+1)}
        +\frac{h}{2}(\rho U(2))^{1-\alpha_{2}}\gamma(\alpha_{2},1)\right), \\
a_{21}&= -\frac{m_{1}B_{\alpha_{1}}^{-1}}{\Gamma(\alpha_{1})}\left(\frac{h^{\alpha_{1}}}{\alpha_{1}(\alpha_{1}+1)}
        +\frac{h}{2}(\rho U(1))^{1-\alpha_{1}}\gamma(\alpha_{1},1)\right),\\
a_{22}&=\frac{m_{2}B_{\alpha_{2}}^{-1}}{\Gamma(\alpha_{2})}\left(\frac{h^{\alpha_{2}}}{\alpha_{2}(\alpha_{2}+1)}
        +\frac{h}{2}(\rho U(2))^{1-\alpha_{2}}\gamma(\alpha_{2},1)\right)- \frac{\rho U(2)h}{2},
\end{align*}
and
\begin{align*}
b_{1}^{\langle0\rangle}=
          &\frac{m_{1}B_{\alpha_{1}}^{-1}}{\Gamma(\alpha_{1})}\left(\frac{h^{\alpha_{1}}}{\alpha_{1}(\alpha_{1}+1)}\alpha_{1}e^{-\rho U(1)h}
            +\frac{h}{2}(\rho U(1))^{1-\alpha_{1}}\gamma(\alpha_{1},1)\right)G_{1,0}-\frac{\rho U(1)h-2}{2}G_{1,0}\\
          & -\frac{m_{2}B_{\alpha_{2}}^{-1}}{\Gamma(\alpha_{2})}\left(\frac{h^{\alpha_{2}}}{\alpha_{2}(\alpha_{2}+1)}\alpha_{2}e^{-\rho U(2)h}
            +\frac{h}{2}(\rho U(2))^{1-\alpha_{2}}\gamma(\alpha_{2},1)\right)G_{2,0},
\end{align*}
\begin{align*}
b_{2}^{\langle0\rangle}=
          &\frac{m_{2}B_{\alpha_{2}}^{-1}}{\Gamma(\alpha_{2})}\left(\frac{h^{\alpha_{2}}}{\alpha_{2}(\alpha_{2}+1)}\alpha_{2}e^{-\rho U(2)h}
            +\frac{h}{2}(\rho U(2))^{1-\alpha_{2}}\gamma(\alpha_{2},1)\right)G_{2,0}-\frac{\rho U(2)h-2}{2}G_{2,0}\\
          &- \frac{m_{1}B_{\alpha_{1}}^{-1}}{\Gamma(\alpha_{1})}\left(\frac{h^{\alpha_{1}}}{\alpha_{1}(\alpha_{1}+1)}\alpha_{1}e^{-\rho U(1)h}
            +\frac{h}{2}(\rho U(1))^{1-\alpha_{1}}\gamma(\alpha_{1},1)\right)G_{1,0}.
\end{align*}
Then, there holds
\begin{displaymath}
\begin{bmatrix}1&0\\0&1\end{bmatrix}\begin{bmatrix}G_{1}(t_{1})\\G_{2}(t_{1})\end{bmatrix}
=\begin{bmatrix}b_{1}^{\langle0\rangle}\\b_{2}^{\langle0\rangle}\end{bmatrix}
+\begin{bmatrix}a_{11}&a_{12}\\a_{21}&a_{22}\end{bmatrix}
\cdot\begin{bmatrix}G_{1}(t_{1})\\G_{2}(t_{1})\end{bmatrix}.
\end{displaymath}
Furthermore, we can obtain that
\begin{equation}\label{for:n0}
\begin{bmatrix}G_{1}(t_{1})\\G_{2}(t_{1})\end{bmatrix}
=\begin{bmatrix}1-a_{11}&-a_{12}\\-a_{21}&1- a_{22}\end{bmatrix}^{-1}\cdot\begin{bmatrix}b_{1}^{\langle0\rangle}\\b_{2}^{\langle0\rangle}\end{bmatrix}.
\end{equation}
Since the matrix $\begin{bmatrix}1-a_{11}&-a_{12}\\-a_{21}&1- a_{22}\end{bmatrix}$
is a principally diagonally dominant matrix, so it is invertible.

\textbf{Case II:} For $n\geq1$, there are
 \begin{align*}
G_{1}(t_{n+1})=&
      \frac{m_{1}B_{\alpha_{1}}^{-1}}{\Gamma(\alpha_{1})}\left(\frac{h^{\alpha_{1}}}{\alpha_{1}(\alpha_{1}+1)}\sum\limits_{j=0}^{n}d_{j,n}^{\langle 1\rangle}G_{1}(t_{j})+(\rho U(1))^{1-\alpha_{1}}\gamma(\alpha_{1},1)\left(\frac{h}{2}G_{1,0}+h\sum\limits_{j=1}^{n}G_{1}(t_{j})\right)\right)\\
      &+\frac{m_{1}B_{\alpha_{1}}^{-1}}{\Gamma(\alpha_{1})}\left(\frac{h^{\alpha_{1}}}{\alpha_{1}(\alpha_{1}+1)}
        +\frac{h}{2}(\rho U(1))^{1-\alpha_{1}}\gamma(\alpha_{1},1)\right)G_{1}(t_{n+1})\\
      &-\frac{m_{2}B_{\alpha_{2}}^{-1}}{\Gamma(\alpha_{2})}\left(\frac{h^{\alpha_{2}}}{\alpha_{2}(\alpha_{2}+1)}\sum\limits_{j=0}^{n}d_{j,n}^{\langle
        2\rangle}G_{2}(t_{j})+(\rho U(2))^{1-\alpha_{2}}\gamma(\alpha_{2},1)\left(\frac{h}{2}G_{2,0}+h\sum\limits_{j=1}^{n}G_{2}(t_{j})\right)\right)\\
      &-\frac{m_{2}B_{\alpha_{2}}^{-1}}{\Gamma(\alpha_{2})}\left(\frac{h^{\alpha_{2}}}{\alpha_{2}(\alpha_{2}+1)}
        +\frac{h}{2}(\rho U(2))^{1-\alpha_{2}}\gamma(\alpha_{2},1)\right)G_{2}(t_{n+1})\\
      &-\rho U(1)\left(\frac{h}{2}G_{1,0}+h\sum\limits_{j=1}^{n}G_{1}(t_{j})\right)-\frac{\rho U(1)h}{2}G_{1}(t_{n+1})+G_{1,0}
\end{align*}
and
 \begin{align*}
G_{2}(t_{n+1})=
      &\frac{m_{2}B_{\alpha_{2}}^{-1}}{\Gamma(\alpha_{2})}\left(\frac{h^{\alpha_{2}}}{\alpha_{2}(\alpha_{2}+1)}\sum\limits_{j=0}^{n}d_{j,n}^{\langle
        2\rangle}G_{2}(t_{j})+(\rho U(2))^{1-\alpha_{2}}\gamma(\alpha_{2},1)\left(\frac{h}{2}G_{2,0}+h\sum\limits_{j=1}^{n}G_{2}(t_{j})\right)\right)\\
      &+\frac{m_{2}B_{\alpha_{2}}^{-1}}{\Gamma(\alpha_{2})}\left(\frac{h^{\alpha_{2}}}{\alpha_{2}(\alpha_{2}+1)}
        +\frac{h}{2}(\rho U(2))^{1-\alpha_{2}}\gamma(\alpha_{2},1)\right)G_{2}(t_{n+1})\\
      &-\frac{m_{1}B_{\alpha_{1}}^{-1}}{\Gamma(\alpha_{1})}\left(\frac{h^{\alpha_{1}}}{\alpha_{1}(\alpha_{1}+1)}\sum\limits_{j=0}^{n}d_{j,n}^{\langle
        1\rangle}G_{1}(t_{j})+(\rho U(1))^{1-\alpha_{1}}\gamma(\alpha_{1},1)\left(\frac{h}{2}G_{1,0}+h\sum\limits_{j=1}^{n}G_{1}(t_{j})\right)\right)\\
      &-\frac{m_{1}B_{\alpha_{1}}^{-1}}{\Gamma(\alpha_{1})}\left(\frac{h^{\alpha_{1}}}{\alpha_{1}(\alpha_{1}+1)}
        +\frac{h}{2}(\rho U(1))^{1-\alpha_{1}}\gamma(\alpha_{1},1)\right)G_{1}(t_{n+1})\\
      &-\rho U(1)\left(\frac{h}{2}G_{2,0}+h\sum\limits_{j=1}^{n}G_{2}(t_{j})\right)-\frac{\rho U(1)h}{2}G_{2}(t_{n+1})+G_{2,0},
\end{align*}
where
\begin{equation}\label{eq:PCa}
d_{j,n}^{\langle m\rangle}=
\left\{
 \begin{aligned}
  &e^{-\rho U(m)(n+1)h}\left(n^{\alpha_{m}+1}-(n+1)^{\alpha_{m}}(n-\alpha_{m})\right),        & if\ j&=0,\\
  &e^{-\rho U(m)(n-j+1)h}\left((n-j+2)^{\alpha_{m}+1}+(n-j)^{\alpha_{m}+1}-2(n-j+1)^{\alpha_{m}+1}\right),   & if\ 1&\leq j\leq n.\\
 \end{aligned}\right.,~~ m=1,2.
\end{equation}
Denote
\begin{displaymath}\begin{aligned}
b_{1}^{\langle n\rangle}=&
\frac{m_{1}B_{\alpha_{1}}^{-1}}{\Gamma(\alpha_{1})}\left(\frac{h^{\alpha_{1}}}{\alpha_{1}(\alpha_{1}+1)}\sum\limits_{j=0}^{n}d_{j,n}^{\langle1\rangle}
    G_{1}(t_{j})+(\rho U(1))^{1-\alpha_{1}}\gamma(\alpha_{1},1)\left(\frac{h}{2}G_{1,0}+h\sum\limits_{j=1}^{n}G_{1}(t_{j})\right)\right)\\
  &-\frac{m_{2}B_{\alpha_{2}}^{-1}}{\Gamma(\alpha_{2})}\left(\frac{h^{\alpha_{2}}}{\alpha_{2}(\alpha_{2}+1)}\sum\limits_{j=0}^{n}d_{j,n}^{\langle2\rangle}
    G_{2}(t_{j})+(\rho U(2))^{1-\alpha_{2}}\gamma(\alpha_{2},1)\left(\frac{h}{2}G_{2,0}+h\sum\limits_{j=1}^{n}G_{2}(t_{j})\right)\right)\\
  &-\rho U(1)\left(\frac{h}{2}G_{1,0}+h\sum\limits_{j=1}^{n}G_{1}(t_{j})\right)+G_{1,0},
\end{aligned}\end{displaymath}
and
\begin{displaymath}\begin{aligned}
b_{2}^{\langle n\rangle}=
  &\frac{m_{2}B_{\alpha_{2}}^{-1}}{\Gamma(\alpha_{2})}\left(\frac{h^{\alpha_{2}}}{\alpha_{2}(\alpha_{2}+1)}\sum\limits_{j=0}^{n}d_{j,n}^{\langle2\rangle}
    G_{2}(t_{j})+(\rho U(2))^{1-\alpha_{2}}\gamma(\alpha_{2},1)\left(\frac{h}{2}G_{2,0}+h\sum\limits_{j=1}^{n}G_{2}(t_{j})\right)\right)\\
 &-\frac{m_{1}B_{\alpha_{1}}^{-1}}{\Gamma(\alpha_{1})}\left(\frac{h^{\alpha_{1}}}{\alpha_{1}(\alpha_{1}+1)}\sum\limits_{j=0}^{n}d_{j,n}^{\langle1\rangle}
    G_{1}(t_{j})+(\rho U(1))^{1-\alpha_{1}}\gamma(\alpha_{1},1)\left(\frac{h}{2}G_{1,0}+h\sum\limits_{j=1}^{n}G_{1}(t_{j})\right)\right)\\
 &-\rho U(2)\left(\frac{h}{2}G_{2,0}+h\sum\limits_{j=1}^{n}G_{2}(t_{j})\right)+G_{2,0}.
\end{aligned}\end{displaymath}
Thus,
\begin{equation}\label{for:nn}
\begin{bmatrix}G_{1}(t_{n+1})\\G_{2}(t_{n+1})\end{bmatrix}
=\begin{bmatrix}1-a_{11}&-a_{12}\\-a_{21}&1- a_{22}\end{bmatrix}^{-1}\cdot\begin{bmatrix}b_{1}^{\langle n\rangle}\\b_{2}^{\langle n\rangle}
\end{bmatrix},~ ~n=1,2,\cdots,M.
\end{equation}

Combining (\ref{for:n0}) with (\ref{for:nn}), the time marching scheme for system (\ref{eq:DIVFeynman-kac}) is obtained.

\section*{Acknowledgments}

The authors Ma and Deng are supported by the National Natural Science Foundation of China under Grant No. 12071195, the AI and Big Data Funds under Grant No. 2019620005000775, and the Innovative Groups of Basic Research in Gansu Province under Grant No. 22JR5RA391. The author Zhao is supported by Guangdong Basic and Applied Basic Research Foundation under Grant No. 2022A1515011332. The authors have no relevant financial or non-financial interests to disclose.

\section*{References}


\begin{thebibliography}{99}


\bibitem[1]{Talbot79} A. Talbot, The accurate numerical inversion of Laplace transforms, IMA J. Appl. Math., 23 (1979), 97-120.

\bibitem[2]{Li22} B. Y. Li and S. Ma, Exponential convolution quadrature for nonlinear subdiffusion equations with nonsmooth initial data, SIAM. J. Numer. Anal., 60  (2022), 503-528.

\bibitem[3]{Li21} B. Y. Li and S. Ma, A high-order exponential integrator for nonlinear parabolic equations with nonsmooth initial data, J. Sci. Comput., 87 (2021), pp. 1-16.

\bibitem[4]{Li23} B. Y. Li, Y. P. Lin, S. Ma, Q. Q Rao, An exponential spectral method using VP means for semilinear subdiffusion equations with rough data, SIAM J. Numer. Anal., 2023 (accepted).

\bibitem[5]{Dingfelder15} B. Dingfelder and J. A. C. Weideman., An improved Talbot method for numerical Laplace transform iverson, Numer. Algorithms, 68 (2015), 167-183.

\bibitem[6]{SIAM18} B. T. Jin, B. Y. Li, and Z. Zhou, Numerical analysis of nonlinear subdiffusion equations, SIAM J. Numer. Anal., 56 (2018), 1-23.

\bibitem[7]{Sloan99} D. Sheen, I. H. Sloan, and V. Thom¨¦e, A parallel method for time-discretization of parabolic problems based on contour integral representation and quadrature, Math. Comp., 69 (2000), 177-195.

\bibitem[8]{Ma22} F.~G. Ma, L.~J. Zhao, W.~H. Deng, and Y.~J. Wang, Analyses of the contour integral method for time fractional subdiffusion-normal transport equation. \href{https://arxiv.org/abs/2210.09594}{arXiv:2210.09594v1}.

\bibitem[9]{Zeng22} H. Zhang, F. H. Zeng, X. Y. Jiang, and G. E. Karniadakis, Convergence analysis of the time-stepping numerical methods for time-fractional nonlinear subdiffusion equations, Fract. Calc. Appl. Anal., 25 (2022), 453-487.

\bibitem[10]{Wang95} H. Wang, R. E. Ewing, and T. F. Russell, Eulerian-Lagrangian localized adjoint methods for convection-diffusion equations and their convergence analysis, IMA J. Numer. Anal., 15 (1995), 405-459.

\bibitem[11]{Podlubny99} I. Podlubny, Fractional Differential Equations, Academic Press, San Diego, 1999.

\bibitem[12]{Lee06}  J. Lee and D. Sheen, A parallel method for backward parabolic problems based on the Laplace transformation, SIAM J. Numer. Anal., 44 (2006), 1466-1486.

\bibitem[13]{Weidemantrefe07} J. A. C. Weideman and L. N. Trefethen, Parabolic and hyperbolic contours for computing the Bromwich integral, Math. Comp., 76 (2007), 1341-1356.

\bibitem[14]{Freed02}  K. Diethelm, N. J. Ford, and A. D. Freed, A predictor-corrector approach for the numerical solution of fractional differential equations, Nonlinear Dynam., 29 (2002), 3-22.

\bibitem[15]{int11} K. J. in¡¯t Hout and J. A. C. Weideman, A contour integral method for the Black-Scholes and Heston equations, SIAM J. Sci. Comput., 33 (2011), 763-785.

\bibitem[16]{Trefethren} L. N. Trefethen and J. A. C. Weideman, The exponentially convergent trapezoidal rule, SIAM Rev., 56 (2014), 385-458.

\bibitem[17]{zhao21} L. J. Zhao, W. H. Deng, and J. S. Hethaven, Characterization of image spaces of riemann-liouville fractional integral operators on sobolev spaces $W^{m,p}(\Omega)$, Sci. China Math., 64 (2021), 2611-2636.

\bibitem[18]{Fernandez04} M. L. Fern\'{a}ndz and C. Palencia, On the numerical inversion of the Laplace transform of certain holomorphic mappings, Appl. Numer. Math., 51 (2004), 289-303.

\bibitem[19]{Fernandez06} M. L. Fern\'{a}ndz, C. Palencia, and A. Sch\"{a}dle, A spectral order method for inverting sectorial Laplace transforms, SIAM J. Numer. Anal., 44 (2006), 1332-1350.

\bibitem[20]{Taiwo95} O. Taiwo, J. Schultz, and V. Krebs, A comparison of two methods for the numerical inversion of Laplace transforms, Comput. Chem. Eng., 19 (1995), 303-308.

\bibitem[21]{xu18} P. B. Xu and W. H. Deng, Fractional compound Poisson process with multiple internal states, Math. Model. Nat. Phenom., 13 (2018), 10.

\bibitem[22]{Piessens75} R. Piessens, A bibliography on numerical inversion of the Laplace transform and applications, J. Comput. Appl. Math., 1 (1975), 115-128.

\bibitem[23]{Piessens76} R. Piessens and N. D. P. Dang, A bibliography on numerical inversion of the Laplace transform and applications: A supplement, J. Comput. Appl. Math., 2 (1976), 225-228.

\bibitem[24]{SLE_07}  S. Carmi, L. Turgeman, and E. Barkai, On distributions of functionals of anomalous diffusion paths, J. Stat. Phys., 141 (2010), 1071-1092.

\bibitem[25]{Deng07JCAM} W. H. Deng, Short memory principle and a predictor-corrector approach for fractional differential equations, J. Comput. Appl. Math., 206 (2007), 174-188.

\bibitem[26]{Dengli18} W. H. Deng, B. Y. Li, Z. Qian, and H. Wang, Time discretization of a tempered fractional Feynman-Kac equation with measure data, SIAM J. Numer. Anal., 56 (2018), 3249-3275.

\bibitem[27]{Sheen06} W. McLean, I. H. Sloan, and V. Thom¨¦e, Time discretization via Laplace transformation of an integro-differential equation of parabolic type, Numer. Math., 102 (2006), 497-522.

\bibitem[28]{Mclean10} W. Mclean and V. Thom¨¦e, Maximum-norm error analysis of a numerical solution via Laplace transformation and quadrature of a fractional order evolution equation, IMA J. Numer. Anal., 30 (2010), 208-230.

\bibitem[29]{wudeng16} X. C. Wu, W. H. Deng, and E. Barkai, Tempered fractional Feynman-Kac equation: theory and examples, Phys. Rev. E, 93 (2016), 032151.

\bibitem[30]{XiYang18} X. Yang, Numerical contour integral methods for unsteady Stokes equations, Comput. Math. Appl., 75 (2018), 4414-4426.

\bibitem[31]{ZhouSun18} Z. Q. Zhou, J. T. Ma, and H. W. Sun, Fast Laplace transform methods for free-boundary problems of fractional diffusion equations, J. Sci. Comput., 74 (2018), 49-69.


















\end{thebibliography}
\end{document}